\documentclass[a4paper,12pt]{amsart}
\usepackage{amssymb}
\usepackage[colorlinks=true]{hyperref}
\usepackage{pstricks}
\usepackage{graphicx}
\usepackage{mathdots}
\allowdisplaybreaks[1]
\newtheorem{lemma}{Lemma}
\newtheorem{proposition}{Proposition}
\newtheorem{theorem}{Theorem}
\newtheorem*{theorem*}{Theorem}
\newtheorem*{corollary*}{Corollary}
\newcommand{\ASM}{\mathrm{ASM}}
\newcommand{\DPP}{\mathrm{DPP}}
\newcommand{\SVDWBC}{\mathrm{6VDW}}
\newcommand{\NILP}{\mathrm{NILP}}
\newcommand{\OT}{\mathrm{OT}}
\newcommand{\paths}{\mathcal{P}}
\newcommand{\ba}[1]{\begin{array}{@{}#1@{}}}
\newcommand{\ea}{\end{array}}
\newcommand{\ds}{\displaystyle}
\newcommand{\ts}{\textstyle}
\renewcommand{\ss}{\scriptstyle}
\newcommand{\sss}{\scriptscriptstyle}
\newcommand{\pict}[2]{\makebox(0,0)[#1]{$#2$}}
\newcommand{\pictss}[2]{\makebox(0,0)[#1]{$\ss#2$}}

\newcommand{\vertex}{\psline[linewidth=0.5pt](0,-0.5)(0,0.5)\psline[linewidth=0.5pt](-0.5,0)(0.5,0)}
\author[R.~E.~Behrend]{Roger E.~Behrend}
\address{R.~E.~Behrend, School of Mathematics, Cardiff University, Cardiff, CF24 4AG, UK}
\email{behrendr@cardiff.ac.uk}
\gdef\s{s}
\title[On the weighted enumeration of ASM{\s} and DPP{\s}]{On the weighted enumeration of
alternating sign matrices and descending plane partitions}
\author[P.~Di Francesco]{Philippe Di Francesco}
\address{P.~Di Francesco,
Institut de Physique Th\'eorique de Saclay,
CEA/DSM/SPhT, CNRS URA 2306,
C.E.A.-Saclay, F-91191 Gif sur Yvette Cedex, France}
\email{philippe.di-francesco@cea.fr}
\author[P.~Zinn-Justin]{Paul Zinn-Justin}
\address{P.~Zinn-Justin, UPMC Univ Paris 6, CNRS UMR 7589, LPTHE,
75252 Paris Cedex, France}
\email{pzinn@lpthe.jussieu.fr}
\thanks{PDF and PZJ acknowledge partial support
from ANR program ``GRANMA'' BLAN08-1-13695.}
\begin{document}
\begin{abstract}
We prove a conjecture of Mills, Robbins and Rumsey [Alternating sign matrices and descending plane partitions,
J.~Combin.\ Theory Ser.~A \textbf{34} (1983), 340--359]
that, for any $n$, $k$, $m$ and $p$,
the number of $n\times n$ alternating sign matrices (ASMs) for which the~1 of the first row is in column~$k+1$ and
there are exactly $m$ $-1$'s and $m+p$ inversions is equal to the number of descending plane partitions (DPPs)
for which each part is at most~$n$ and there are exactly $k$ parts equal to~$n$,~$m$ special parts and~$p$ nonspecial parts.
The proof involves expressing the associated generating functions for ASMs and DPPs with fixed $n$ as
determinants of $n\times n$ matrices, and using elementary transformations to show that these determinants are equal.
The determinants themselves are obtained by standard methods: for ASMs this involves
using the Izergin--Korepin formula for the partition function of the six-vertex model with domain-wall boundary conditions, together
with a bijection between ASMs and configurations of this model, and for DPPs it involves using
the Lindstr\"{o}m--Gessel--Viennot theorem, together with a bijection between
DPPs and certain sets of nonintersecting lattice paths.
\end{abstract}
\maketitle

\section{Preliminaries}
Alternating sign matrices (ASMs) and descending plane partitions (DPPs) are combinatorial objects
which arose within a few years of each other in the late 1970s and early 1980s, but in somewhat different contexts.
DPPs were introduced by
Andrews while attempting to prove a conjectured formula for the generating function of cyclically symmetric plane
partitions~\cite{And-weakMac,And-MacDPP},
whereas ASMs first appeared during studies by Mills, Robbins and Rumsey of Dodgson's condensation algorithm for the evaluation of
determinants~\cite{MRR-Mac,MRR-ASM,RR-ASM}.

However, despite their independent origins, it soon became apparent
that ASMs and DPPs share some basic enumerative properties.  In
particular, it was conjectured by Mills, Robbins and
Rumsey~\cite[Conj.~3]{MRR-ASM}, in one of the early papers on the subject, that certain finite sets of ASMs have
the same sizes as certain finite sets of DPPs, where these sets are
comprised of all ASMs or DPPs with fixed values of particular
statistics. It is the objective of this paper to prove this conjecture.

Some special cases of this result are already known to be valid, the most noteworthy perhaps
being the fact, which follows from results of Andrews~\cite{And-MacDPP} and of Zeilberger~\cite{Zeil-ASM} or Kuperberg~\cite{Kup-ASM}, that
for any positive integer~$n$, the number of $n\times n$ ASMs is equal to the number of DPPs in which each part is at most~$n$,
with these numbers furthermore being given by a simple product formula.

Many other similar results, involving particular classes, refinements, representations or weightings of ASMs or plane partitions,
and equalities of sizes of finite sets or explicit formulae for these sizes, have also been proved.
Such results can typically be stated very simply, and yet are often conjectured long before they are proved.

Ideally, results of this sort would be proved bijectively, but unfortunately few such proofs are currently known in this area.
Nevertheless, there are some alternative methods of proof which have been widely successful.
For ASMs, one of the most successful techniques involves first applying a simple correspondence
between ASMs and configurations of certain cases of the statistical mechanical six-vertex model, then
using the integrability of this model to obtain a determinant or Pfaffian expression for the model's partition function,
and finally using certain methods related to determinants or Pfaffians to transform or evaluate the determinant or Pfaffian,
for particular values of the additional parameters which necessarily appear, thereby leading to
enumerative results.
For plane partitions, one of the most successful techniques involves first applying a simple correspondence
between plane partitions and sets of nonintersecting lattice paths, then using certain general results for nonintersecting paths
to express the number of such sets of paths,
or the associated generating function, as a determinant or Pfaffian, and finally again
using certain methods to transform or evaluate the determinant or Pfaffian,
thus leading to enumerative results.

In this paper, elements of each of these by-now-standard techniques are employed to obtain determinant formulae for the
generating functions associated with the sets of ASMs and DPPs under consideration.
Transformations are then applied to the generating functions associated with the entries of the underlying matrices
of these determinants, and this enables it to be shown that the determinants are equal,
which in turn implies that the sizes of the sets in the Mills, Robbins and Rumsey conjecture are equal, but
without these sizes having been evaluated explicitly.  Although this proof has much in common
with other known proofs in this area, one of its distinguishing features
is that it involves deriving a determinant for the partition function of
the relevant case of the six-vertex model
which depends directly on the underlying weights of the model,
rather than on the spectral parameters of the model (which have only an implicit inverse-functional dependence
on the weights, resulting from certain rational functions which give the weights in terms of the spectral parameters).
Indeed, once this determinant for the partition function
in terms of the weights has been obtained, the equality of the determinants for ASMs and DPPs
follows relatively straightforwardly.

For reviews of, and historical information about, enumerative results and conjectures involving
ASMs or plane partitions, see for example Bressoud~\cite{Bressoud}, Bressoud and Propp~\cite{BP-ASM},
Propp~\cite{Propp-match,Propp-ASM,Propp-matchupdate}, Robbins~\cite{Rob-story,Rob-ASM},
Stanley \cite{Stan,Stan-symPP,Stan-symPPerr} and Zeilberger~\cite{Zeil-Robbins}.
Note, however, that most of these reviews were written over a decade ago,
and that many of the conjectures which they report have since been proved.
For example, most of the conjectures for symmetry classes of ASMs listed by Robbins in~\cite[Table~3]{Rob-story} have been proved by
Kuperberg~\cite{Kup-symASM}, Okada~\cite{Oka} or Razumov and Stroganov~\cite{RS-halfturn,RS-quarterturn}.
Note also that, in the past decade, further interesting results and conjectures
involving ASMs or plane partitions,
and revealing new connections with algebra, mathematical physics and other areas,
have appeared, the most prominent of these being the Razumov--Stroganov (ex-)conjecture.
For information about, and references for, these more recent developments see, for example, the reviews
of de Gier~\cite{dG-review,dG-review2,dG-review3} and Zinn-Justin~\cite{hdr}, and the proof of the
Razumov--Stroganov conjecture by Cantini and Sportiello~\cite{CS-RS}.

An outline of this paper is as follows.  In Sections~\ref{defin} and~\ref{stat} we give the definitions of ASMs and DPPs,
of certain statistics for ASMs and DPPs, and of the generating
functions associated with these statistics.  In Section~\ref{MRRConj} we state the Mills, Robbins and Rumsey conjecture,
and outline the cases which have been previously proved or studied.
We then present the proof of this conjecture in the next two sections.
For ease of understanding, the proof of a certain, unrefined version of the result is
given first in Section~\ref{unref}, with the proof of the
general, refined case following in Section~\ref{refenum}.  Both of these sections have
the same structure, in which the ASM generating function is expressed as a determinant
in Sections~\ref{ASMdetderiv} and~\ref{ASMdetderivref},
the DPP generating function is expressed as a determinant in Sections~\ref{DPPdetderiv} and~\ref{DPPdetderivref},
and it is shown that the determinants are equal in Sections~\ref{ASMdetDPPdetproof} and~\ref{ASMdetDPPdetproofref}.
We conclude the paper by discussing some further aspects of this work, and some future work, in Section~\ref{discuss}.

\subsection{Definitions}\label{defin}
In this subsection, definitions will be given for ASMs and DPPs, and for sets of ASMs and DPPs of order $n$.

An ASM, first defined by Mills, Robbins and Rumsey~\cite{MRR-Mac,MRR-ASM},
is a square matrix in which:
\begin{list}{$\bullet$}{\setlength{\labelwidth}{2mm}\setlength{\leftmargin}{10mm}}
\item Each entry is~$0$, $1$ or~$-1$.
\item The nonzero entries alternate in sign along each row and column.
\item The sum of entries in each row and column is $1$.
\end{list}
It follows that an ASM has a unique 1 in each of its first and last row and column,
and that any permutation matrix is an ASM.
A further example of an ASM is
\begin{equation}\label{ASMEx}A=\begin{pmatrix}
0&0&0&1&0&0\\0&1&0&-1&1&0\\1&-1&1&0&0&0\\0&0&0&1&0&0\\0&1&0&-1&0&1\\0&0&0&1&0&0\end{pmatrix}.\end{equation}
A DPP, first defined by Andrews~\cite{And-weakMac,And-MacDPP}, is an array of positive integers, called parts,
of the form
\setlength{\unitlength}{5pt}
\begin{equation}\label{DPP}\ba{c@{\:}c@{\:}c@{\:}l}
D_{11}&D_{12}&D_{13}&\begin{picture}(32.1,0)\multiput(0,0)(1,0){32}{.}\end{picture}D_{1,\lambda_1}\\
&D_{22}&D_{23}&\begin{picture}(27.1,0)\multiput(0,0)(1,0){27}{.}\end{picture}D_{2,\lambda_2+1}\\
&&D_{33}&\begin{picture}(22.1,0)\multiput(0,0)(1,0){22}{.}\end{picture}D_{3,\lambda_3+2}\\
&&&\begin{picture}(20,3.3)\multiput(1,0)(-1,1){3}{.}\multiput(19,0)(1,1){3}{.}\end{picture}\\
&&&\hspace*{8pt}D_{tt}\begin{picture}(8.6,0)\multiput(0.5,0)(1,0){8}{.}\end{picture}D_{t,\lambda_t+t-1}\ea\end{equation}
in which:
\begin{list}{$\bullet$}{\setlength{\labelwidth}{2mm}\setlength{\leftmargin}{5mm}}
\item The parts decrease weakly along rows, i.e., $D_{ij}\ge D_{i,j+1}$ whenever both sides are defined.
\item The parts decrease strictly down columns, i.e., $D_{ij}>D_{i+1,j}$ whenever both sides are defined.
\item The first parts of each row and the row lengths satisfy
\begin{equation}\hspace*{5mm}D_{11}>\lambda_1\ge D_{22}>\lambda_2\ge\ldots\ge D_{t-1,t-1}>\lambda_{t-1}\ge D_{tt}>\lambda_t.\end{equation}
\end{list}
The empty array is regarded as a DPP and denoted $\emptyset$.
A further example of a DPP is
\begin{equation}\label{DPPEx}\ba{c@{\;}c@{\;}c@{\;}c@{\;}c@{\;}c}
&6&6&6&5&2\\
D=&&4&4&1\\
&&&3\ea.\end{equation}
We now define sets of ASMs and DPPs of order $n$, for each positive integer $n$.
Let $\ASM(n)$ be the set of all $n\times n$ ASMs, and let $\DPP(n)$ be the set of all
DPPs in which each part is at most~$n$.

\makeatletter\newcommand{\biggg}{\bBigg@{3}}\makeatother
For example,
\begin{alignat}{2}\label{ASM3}
\ASM(3)&=&\biggg\{&\begin{pmatrix}1&0&0\\0&1&0\\0&0&1\end{pmatrix},
\begin{pmatrix}0&0&1\\0&1&0\\1&0&0\end{pmatrix},
\begin{pmatrix}1&0&0\\0&0&1\\0&1&0\end{pmatrix},
\begin{pmatrix}0&0&1\\1&0&0\\0&1&0\end{pmatrix},\\[1.5mm]
&&&\begin{pmatrix}0&1&0\\1&0&0\\0&0&1\end{pmatrix},
\begin{pmatrix}0&1&0\\0&0&1\\1&0&0\end{pmatrix},
\begin{pmatrix}0&1&0\\1&-1&1\\0&1&0\end{pmatrix}\biggg\},\notag\\[3mm]
\label{DPP3}\DPP(3)&=\rlap{$\left\{\emptyset,\;\ba{c@{\;}c}3&3\\[-0.8mm]&2\ea,\;2,\;3\;3,\;3,
\;3\;2,\;3\;1\right\}$.}\end{alignat}
Note that a DPP~$D$
will often be associated with a particular value of~$n$ for which $D\in\DPP(n)$ (i.e., a particular $n\ge D_{11}$).
For example, the DPP of~\eqref{DPPEx}
will be associated throughout this paper with $n=6$.

\subsection{Statistics and generating functions}\label{stat}
In this subsection, certain statistics will be introduced for ASMs and DPPs,
and the associated generating functions will be defined.

For a given positive integer $n$, define statistics for each $A\in\ASM(n)$ as
\begin{align}
\ds\label{nuA}\nu(A)&=\sum_{\substack{1\le i<i'\le n\\1\le j'\le j\le n}}\!A_{ij}\,A_{i'j'}=
\sum_{\substack{1\le i\le i'\le n\\1\le j'<j\le n}}\!A_{ij}\,A_{i'j'},\\[1mm]
\label{muA}\mu(A)&=\text{number of $-1$'s in $A$},\\[1mm]
\label{rhoA}\rho(A)&=\text{number of 0's to the left of the 1 in the first row of $A$,}\\
\intertext{and define statistics for each $D\in\DPP(n)$ as}
\ds\label{nuD}\nu(D)&=\text{number of parts of $D$ for which $D_{ij}>j-i$},\\[1mm]
\label{muD}\mu(D)&=\text{number of parts of $D$ for which $D_{ij}\le j-i$},\\[1mm]
\label{rhoD}\rho(D)&=\text{number of parts equal to $n$ in (necessarily the first row of) $D$}.
\end{align}
The equality between the two expressions for $\nu(A)$ in~\eqref{nuA} holds
for any matrix $A$ in which the sum of entries in a row and column is the same for all rows and columns,
and therefore holds for an ASM $A$.

A part of a DPP $D$ for which $D_{ij}\le j-i$ is referred to by Mills, Robbins and Rumsey~\cite[p.~344]{MRR-ASM} as
a special part.  Thus, $\nu(D)$ and $\mu(D)$ are the numbers of nonspecial and special parts respectively in~$D$.

The examples~\eqref{ASMEx} and~\eqref{DPPEx} have $\nu(A)=5$, $\mu(A)=3$, $\rho(A)=3$,
$\nu(D)=7$, $\mu(D)=2$ and (taking $n=6$) $\rho(D)=3$.

Definitions~\eqref{nuA}--\eqref{rhoD} are based on definitions of statistics for ASMs and DPPs first introduced by
Mills, Robbins and Rumsey~\cite[pp.~344--345]{MRR-ASM}.  However, note that two of the three statistics used in~\cite{MRR-ASM}
differ slightly from those used here.
In particular, Mills, Robbins and Rumsey use statistics $\nu'$, $\mu$ and $\rho'$, where~$\mu$ is again given by~\eqref{muA}
and~\eqref{muD}, while $\nu'$ and $\rho'$ are given by
$\nu'(A)=\sum_{1\le i<i'\le n,\:1\le j'<j\le n}A_{ij}\,A_{i'j'}$ for $A\in\ASM(n)$,
$\nu'(D)=$ total number of parts in $D$ for $D\in\DPP(n)$, and $\rho'(X)=\rho(X)+1$
for $X\in\ASM(n)$ or $X\in\DPP(n)$.  It follows that $\nu'(X)=\nu(X)+\mu(X)$ for any $X\in\ASM(n)$ or $X\in\DPP(n)$.

The statistic $\nu'$ for an ASM $A$ is referred to by Mills, Robbins and Rumsey~\cite[p.~344]{MRR-ASM}
as the number of inversions in~$A$, since it generalizes the usual definition of the number of inversions for permutation matrices.
More specifically, if $A\in\ASM(n)$ is a permutation matrix, i.e., if $\mu(A)=0$, then $\nu'(A)$ is the number of inversions
in the permutation $\pi\in\mathcal{S}_n$ given by $\delta_{\pi_i,j}=A_{ij}$.  In fact, if $A$ is a permutation matrix,
then $\nu'(A)=\nu(A)$ so that $\nu$ could be regarded as an alternative generalization of the number of inversions for permutation matrices.

Now define generating functions, which give weighted enumerations of the elements of $\ASM(n)$ or~$\DPP(n)$ using
arbitrary weights~$x$, $y$ and~$z$ associated
with the statistics \eqref{nuA}--\eqref{rhoA} or~\eqref{nuD}--\eqref{rhoD}, as
\begin{align}\label{ZASM}Z_\ASM(n,x,y,z)&=\sum_{A\in\ASM(n)}x^{\nu(A)}\,y^{\mu(A)}\,z^{\rho(A)},\\
\label{ZDPP}Z_\DPP(n,x,y,z)&=\sum_{D\in\DPP(n)}x^{\nu(D)}\,y^{\mu(D)}\,z^{\rho(D)}.\end{align}
For example, \eqref{ASM3} and~\eqref{DPP3} give
\begin{equation}\label{Z3}Z_\ASM(3,x,y,z)=Z_\DPP(3,x,y,z)=1+x^3z^2+x+x^2z^2+xz+x^2z+xyz,\end{equation}
where the terms are written in an order which corresponds to that used in~\eqref{ASM3} and~\eqref{DPP3}.

It follows easily from the definitions of ASMs and DPPs that the generating functions satisfy
\begin{equation}Z_\ASM(n,x,y,0)=Z_\ASM(n-1,x,y,1),\qquad Z_\DPP(n,x,y,0)=Z_\DPP(n-1,x,y,1).\end{equation}

\subsection{The Mills, Robbins and Rumsey conjecture}\label{MRRConj}
In this subsection, the Mills, Robbins and Rumsey ASM-DPP conjecture will be stated,
and special cases for which the result is already known to be valid will be listed.

The result conjectured by Mills, Robbins and Rumsey~\cite[Conj.~3]{MRR-ASM} (see also Bressoud~\cite[Conj.~10]{Bressoud}),
whose proof is the primary focus of this paper, is as follows.
\begin{theorem}\label{MRRC1}
The sizes of
$\{A\in\ASM(n)\mid\nu(A)=p,\;\mu(A)=m,\;\rho(A)=k\}$
and $\{D\in\DPP(n)\mid\nu(D)=p,\;\mu(D)=m,\;\rho(D)=k\}$ are
equal for any $n$, $p$, $m$ and $k$.

Equivalently,
\begin{equation}\label{MRRC2}Z_\ASM(n,x,y,z)=Z_\DPP(n,x,y,z),\quad
\text{for any $n$, $x$, $y$ and $z$.}\end{equation}
\end{theorem}

Note that the statement of this result in~\cite{MRR-ASM} uses the slightly different statistics $\nu'$ and $\rho'$ outlined in
Section~\ref{stat}.
Note also that the ranges of the integers $p$, $m$ and $k$ can be restricted to
\begin{equation*}\ts p=0,\ldots,\frac{n(n-1)}{2},\qquad m=0,\ldots,\begin{cases}\frac{(n-1)^2}{4},&n\text{ odd}\\
\frac{n(n-2)}{4},&n\text{ even,}\end{cases}\qquad k=0,\ldots,n-1,\end{equation*}
since each set in Theorem~\ref{MRRC1} is otherwise empty.

Furthermore, using certain symmetry operations on ASMs and DPPs which will be defined in Section~\ref{symm}
(and, in particular, applying~\eqref{symstat}),
it follows that the size of either of the sets in Theorem~\ref{MRRC1} is invariant under the replacement of~$p$ by~$n(n-1)/2-p-m$
and~$k$ by~$n-1-k$.  Equivalently, the generating functions satisfy
$Z_\ASM(n,x,y,z)=x^{n(n-1)/2}z^{n-1}Z_\ASM(n,\frac{1}{x},\frac{y}{x},\frac{1}{z})$
and $Z_\DPP(n,x,y,z)=x^{n(n-1)/2}z^{n-1}Z_\DPP(n,\frac{1}{x},\frac{y}{x},\frac{1}{z})$.

The following two sections are devoted to the proof of Theorem~\ref{MRRC1}.
We first prove the unrefined case $z=1$ in Section~\ref{unref}, and then proceed to the
refined case of arbitrary $z$ in Section~\ref{refenum}.  Both of these proofs have an identical structure, in which
the ASM and DPP generating functions are each expressed as determinants, and it is then shown that
the determinants are equal.  However, since the details for $z=1$ are somewhat simpler,
it seemed pedagogically preferable to consider this special case first.
Although the proof of Section~\ref{unref} is thereby repeated in Section~\ref{refenum},
some of the necessary notation and background material appears only in Section~\ref{unref}.

In the remainder of the current section, we list special cases of Theorem~\ref{MRRC1} which have been previously proved or studied.
However, this information is not needed for the full proofs of Sections~\ref{unref}
and~\ref{refenum}, and so could be omitted by readers who are interested only in those proofs.

The cases of Theorem~\ref{MRRC1} which have already been proved in the literature are as follows.
\begin{list}{$\bullet$}{\setlength{\labelwidth}{2mm}\setlength{\leftmargin}{4mm}}
\item\textit{$n$ small, $p$, $m$ and $k$ arbitrary.} The conjecture can easily be verified by hand for small~$n$,
and by computer for somewhat larger~$n$.  For example, Mills, Robbins and Rumsey~\cite[pp.~345--346]{MRR-ASM}
explicitly considered the case $n=5$, $p=3$, $m=1$ and $k=2$, for which there are~10 ASMs and 10~DPPs, and they
stated that they had checked the conjecture by computer for all cases with $n\le7$.
The validity of the conjecture for $n=3$
can also be seen from~\eqref{Z3}.
Note that for $n\le 3$, each triplet of values of~$p$,~$m$ and~$k$ is associated with at most one ASM or DPP, but that for
each $n\ge4$ there are
cases associated with several ASMs or DPPs.
\item\textit{$n$ arbitrary, $p$, $m$ and $k$ summed over (i.e., $x=y=z=1$).}
In this case, which corresponds to straight enumeration of ASMs or DPPs,
\begin{equation}\label{ASMDPPstr}\hspace{4mm}|\ASM(n)|=|\DPP(n)|=\prod_{i=0}^{n-1}\frac{(3i+1)!}{(n+i)!}.\end{equation}
The product formula for $|\ASM(n)|$
was conjectured by Mills, Robbins and Rumsey~\cite[Conj.~1]{MRR-Mac},~\cite[Conj.~1]{MRR-ASM}, and first
proved by Zeilberger~\cite{Zeil-ASM} and Kuperberg~\cite{Kup-ASM}, using two very different methods.
The product formula for $|\DPP(n)|$ (in a slightly different form) was obtained by Andrews~\cite[Thm.~10]{And-weakMac}.
\item\textit{$n$ arbitrary, $p$ and $m$ summed over, $0\le k\le n-1$ (i.e., $x=y=1$).}
In this case, which corresponds to so-called refined enumeration of ASMs or DPPs,
\begin{multline}\hspace{4mm}\label{ASMDPPref}|\{A\in\ASM(n)\mid\rho(A)=k\}|=|\{D\in\DPP(n)\mid\rho(A)=k\}|=\\
\frac{(n+k-1)!\,(2n-k-2)!}{(2n-2)!\:k!\:(n-k-1)!}\;\prod_{i=0}^{n-2}\!\frac{(3i+1)!}{(n+i-1)!}.\end{multline}
The product formula for refined ASM enumeration was conjectured (as an essentially equivalent result)
by Mills, Robbins and Rumsey~\cite[Conj.~2]{MRR-ASM}, and first proved by Zeilberger~\cite{Zeil-refASM}.
(For alternative proofs of this, and the product formula for $|\ASM(n)|$ in~\eqref{ASMDPPstr},
see for example Colomo and Pronko~\cite{CP-orthog1,CP-orthog2} and Fischer~\cite{Fis-RefASM}.)
The product formula for refined DPP enumeration follows from results of Mills, Robbins and Rumsey~\cite{MRR-Mac}. (For further details
see Bressoud~\cite[Conj.~9 and Sec.~5.3]{Bressoud}.)
\item\textit{$n$, $p$ and $k$ arbitrary, $m=0$ (i.e., permutation matrices, DPPs with no special parts).}
A proof for this case was known to Mills, Robbins and Rumsey~\cite[p.~345]{MRR-ASM},
who stated that the conjecture ``has been proved for the case $m=0$'',
but did not provide further details.  (See also Bressoud~\cite[Exercise~6.1.2]{Bressoud}.) However,
a bijection between $\{A\in\ASM(n)\mid\nu(A)=p,\;\mu(A)=0,\;\rho(A)=k\}$
and $\{D\in\DPP(n)\mid\nu(D)=p,\;\mu(D)=0,\;\rho(D)=k\}$
is outlined briefly by Lalonde~\cite{Lal-ASMDPP} and in more detail
(but using slightly different terminology) by Striker~\cite{Stri-ASMDPP}.
The mapping from~$A$ to~$D$ is as follows. First, let
$\chi_i$ be the number of~1's below and to the left of the~1 in row~$i$ of~$A$,
i.e., $\chi_i=\sum_{i<i'\le n,\:1\le j'<j\le n}A_{ij}\,A_{i'j'}$.
(Note that $(\chi_1,\ldots,\chi_n)\in\{0,1,\ldots,n-1\}\times\ldots\times\{0,1\}\times\{0\}$, so that
this is an inversion table.)  Then create~$D$, in the form of~\eqref{DPP}, using $\chi_i$ $(n+1-i)$'s for each $i=1,\ldots,n$,
by taking these parts in weakly decreasing order and forming the rows of~$D$
successively from top to bottom,
placing as many parts as possible successively from left to right into each row, subject to the condition that the
length of a row does not exceed its rightmost part.
The mapping from~$D$ to~$A$ is as follows. First, let~$\chi_i$ be the number of parts of $D$ equal to~$n+1-i$.
Then form a permutation $\pi\in\mathcal{S}_n$ by successively, for $i=1,\ldots,n$,
removing the $(\chi_i+1)$th largest element from $\{1,\ldots,n\}$ and setting this to be $\pi_i$. Finally,
obtain the entries of $A$ as $A_{ij}=\delta_{\pi_i,j}$.
It follows from this bijection that
\begin{equation}\label{Z0}\hspace{4mm}Z_\ASM(n,x,0,z)=Z_\DPP(n,x,0,z)=[n]_{xz}\,[n-1]_x!,\end{equation}
where $[n]_x=1+x+\ldots+x^{n-1}$ and $[n]_x!=[n]_x[n-1]_x\ldots[1]_x$.
An alternative bijection between $\{A\in\ASM(n)\mid\mu(A)=0,\;\rho(A)=k,\;A_{ij}=\delta_{\pi_i,j},\;\pi\text{ has }t\text{ ascents}\}$ and
$\{D\in\DPP(n)\mid\mu(D)=0,\;\rho(D)=k,\;D\text{ has }t\text{ rows}\}$ is given by Ayyer~\cite{Ayy-DPP}, but
for this bijection~$\nu(A)$ does not equal $\nu(D)$ for all corresponding~$A$ and~$D$.
Both of these bijections can be visualized easily in terms of the set~$\NILP(n)$ of sets of nonintersecting lattice paths which will be defined in
Section~\ref{DPPdetderiv}.  Also, both of these bijections satisfy the property
that if~$A$ corresponds to~$D$, then~$A^\ast$ corresponds to~$D^\ast$, where~$\ast$ is a certain symmetry operation on ASMs and DPPs
which will be defined in Section~\ref{symm}.
\item\textit{$n$, $p$ and $k$ arbitrary, $m=1$ (i.e., ASMs with a single $-1$, DPPs with a single special part)}.  This case
was proved nonbijectively by Lalonde in~\cite{Lal-ASM}, with further
aspects being studied in~\cite{Lal-ASMDPP}.  Briefly, the proof involves showing that
$\sum_{A\in\ASM(n),\,\mu(A)=1,\,\rho(A)=k}x^{\nu(A)}$ and $\sum_{D\in\DPP(n),\,\mu(D)=1,\,\rho(D)=k}x^{\nu(D)}$
are each determined by the same recursion relations and boundary conditions.  The relations for ASMs follow from
a simple combinatorial examination of ASMs with a single~$-1$,
while the relations for DPPs are obtained by applying algebraic transformations
to a determinant formula, which is itself obtained using lattice paths and a Lindstr\"{o}m--Gessel--Viennot-type theorem.
\end{list}
Some further cases which do not seem to have been considered explicitly in the literature, but which can
be proved easily, are as follows.
\begin{list}{$\bullet$}{\setlength{\labelwidth}{2mm}\setlength{\leftmargin}{4mm}}
\item\textit{$n$ arbitrary, $p=k(k+1)/2$, $m=k(n-k-1)$, $0\le k\le n-1$.} In this case, there
is a single ASM $A$ and a single DPP $D$ for each~$n$ and~$k$,
where
\begin{equation}\hspace{4mm}A_{ij}=\begin{cases}(-1)^{i+j+k},&k+2\le i+j\le 2n-k\text{ and }|i-j|\le k\\
0,&\text{otherwise,}\end{cases}\end{equation}
and row $i$ of $D$ is comprised of $k-i+1$ $(n-i+1)$'s (all of which are nonspecial parts) followed by $n-k-1$ $(k-i+1)$'s
(all of which are special parts), for  $i=1,\ldots,k$.
\makeatletter\newcommand{\biG}{\bBigg@{8}}\makeatother
Thus,
\begin{equation}\hspace{4mm}\rule{0ex}{14.5ex}A=\biG(\raisebox{0.7ex}{\,$
\ba{c@{\;}c@{}c@{}c@{\,}c@{\;}c@{}c@{\;}c}
\multicolumn{3}{l}{\smash{\!\!\overbrace{\rule{6.5ex}{0ex}}^k}}&&\multicolumn{4}{l}{\smash{\overbrace{\rule{10.5ex}{0ex}}^{n-k-1}}}\\[-1.5mm]
0&\cdots&0&1&0&\cdots&\cdots&0\\[-2mm]
\vdots&\iddots&1&-1&1&&&\vdots\\[-2.6mm]
0&1&-1&&&\ddots&&\vdots\\
1&-1&&&&&1&0\\
0&1&&&&&-1&1\\[-2mm]
\vdots&&\ddots&&&-1&1&0\\[-2.6mm]
\vdots&&&1&-1&1&\iddots&\vdots\\
0&\cdots&\cdots&0&1&0&\cdots&0\ea$\,}\biG),\qquad
D=\ba{c@{\;}c@{\;}c@{}c@{}c@{\;}c@{}c@{\;}c}
\multicolumn{4}{l}{\smash{\!\!\overbrace{\rule{14.5ex}{0ex}}^k}}&\multicolumn{4}{l}{\smash{\,\,\overbrace{\rule{12.7ex}{0ex}}^{n-k-1}}}\\[-1.5mm]
n&n&\cdots&n&k&\cdots&\cdots&k\\
&n\!-\!1&\cdots&n\!-\!1&k\!-\!1&\cdots&\cdots&k\!-\!1\\[-2mm]
&&\ddots&\vdots&\vdots&&&\vdots\\
&&&n\!-\!k\!+\!1&1&\cdots&\cdots&1\ea
\end{equation}
Note that for $k=0$, $A$ is the $n\times n$ identity matrix and $D=\emptyset$, and these are the only ASM and DPP for
which $\nu(A)=\nu(D)=0$.
\item\textit{$n$, $m$ and $k$ arbitrary, $p=1$ (i.e., ASMs with $m+1$ inversions, DPPs with a single nonspecial part).} In this case,
there are $n-m-2$ ASMs and $n-m-2$ DPPs for $k=0$ and each $m=0,\ldots,n-3$,
there is a single ASM and a single DPP for $k=1$ and each $m=0,\ldots,n-2$ ,
and there are no ASMs or DPPs for other values of~$m$ and~$k$.
The ASMs $A$ are given by
\begin{equation}\hspace{4mm}A_{ij}=\begin{cases}1,&1\le i=j\le t-1,\ \ t+m+2\le i=j\le n,\\
&t+1\le i=j+1\le t+m+1\ \ \text{or}\ \ t+1\le i+1=j\le t+m+1\\
-1,&t+1\le i=j\le t+m\\
0,&\text{otherwise,}\end{cases}\end{equation}
and the DPPs consist of
a single row comprised of a single $n-t+1$ followed by $m$ $1$'s, where $t=1,\ldots,n-m-1$.
The cases $t=2,\ldots,n-m-1$ give $k=0$,
while the case $t=1$ gives $k=1$.
\item\textit{$n$, $m$ and $k$ arbitrary, $p=n(n-1)/2-m-1$.}
The ASMs and DPPs for this case can be obtained from those of the previous case by applying certain symmetry operations
which will be defined in Section~\ref{symm}.
\item\textit{$n$, $m$ and $k$ arbitrary, $p=2$ or $p=n(n-1)/2-m-2$.} These cases can be proved
using an approach similar to that used for the previous two cases, although the details become more complicated.
\end{list}
Finally, some closely related cases which have been studied, but without proofs of the conjecture
for all instances of these cases being obtained, are as follows.
\begin{list}{$\bullet$}{\setlength{\labelwidth}{2mm}\setlength{\leftmargin}{4mm}}
\item\textit{$n$ and $m$ arbitrary, $p$ and $k$ summed over.}
In this case, some properties shared by the associated generating functions $Z_\ASM(n,1,y,1)$
and $Z_\DPP(n,1,y,1)$ have been found.  Furthermore, some formulae for these functions at $y=2$ and $y=3$
have been obtained, this often being referred to as 2- or 3-enumeration.  In fact, some formulae for refined 2- and 3-enumeration (i.e., in
which~$k$ is arbitrary rather than being summed over) are also known.

More specifically, a factorization property for $Z_\ASM(n,1,y,1)$ was conjectured by Mills, Robbins and Rumsey~\cite[Conj.~4]{MRR-ASM},
and proved by Kuperberg~\cite[Thm.~3]{Kup-ASM}. (See also Mills, Robbins and Rumsey~\cite[Conj.~5]{MRR-ASM} and
Kuperberg~\cite[Thm.~4, first two equations]{Kup-symASM}.)  A formula for (refined) 2-enumeration of ASMs was obtained by Mills, Robbins
and Rumsey~\cite[Cor.,~p.~358]{MRR-ASM}, and a formula for 3-enumeration of ASMs was conjectured by Mills, Robbins
and Rumsey~\cite[Conj.~6]{MRR-ASM} and proved (up to a certain correction) by Kuperberg~\cite[Thm.~3]{Kup-ASM}.
Most of the corresponding results for DPPs have been, or can be, obtained from results of Mills, Robbins and
Rumsey~\cite[pp.~50 \&~54]{MRR-symPP} which give a determinant formula for $Z_\DPP(n,1,y,1)$.
For further related work (including some refined cases) and references see for example
Colomo and Pronko~\cite{CP-ref3,CP-orthog1,CP-orthog2},
Elkies, Kuperberg, Larsen and Propp~\cite{EKLP1,EKLP2},
Kuperberg~\cite[p.~837]{Kup-symASM},
Mills, Robbins and Rumsey~\cite[Conj.~5]{MRR-TSSCPP},
Okada~\cite[Thm.~2.4(1)]{Oka},
Robbins~\cite[Sec.~2]{Rob-ASM},
Robbins and Rumsey~\cite{RR-ASM},
and Stroganov~\cite{Strog-3enum}.

This case for ASMs has also been studied combinatorially by Cori, Duchon and Le Gac~\cite{CDLG,LeGac}.
Define an isolated~$1$ in an ASM $A$ to be an entry $A_{ij}=1$ such that all other entries in row~$i$ and column~$j$ of~$A$ are~$0$.
It is then shown in~\cite{CDLG,LeGac} that
$|\{A\in\ASM(n)\mid\mu(A)=m\}|=\sum_{i=1}^{3m}\frac{(n!)^2}{(i!)^2(n-i)!}\,C_{i,m}$, where
$C_{i,m}$ is the number of $i\times i$ ASMs which contain $m$ $-1$'s but no isolated $1$'s.
\item\textit{$n$ and $p$ arbitrary, $m$ and $k$ summed over.}
In this case, enumerative formulae for ASMs and DPPs follow from
results of Behrend~\cite[Sec.~14]{Behr-osc}.
Given a Young (Ferrers) diagram~$\lambda$ and a nonnegative integer~$l$, an oscillating (or up-down) tableau of shape $\lambda$ and length~$l$
is defined to be a sequence of $l+1$ Young diagrams which starts with the empty diagram $\emptyset$, ends with~$\lambda$, and in which
successive diagrams differ by the addition or deletion of a single square, without the movement of any other squares.
(The conventions being used here are that a Young diagram consists of an array of adjacent squares justified along the left and top,
with rows numbered from top to bottom, columns numbered from left to right, and weakly-decreasing row lengths.)
Let $\OT(\lambda,l)$ denote the set of all oscillating tableaux of shape~$\lambda$ and length~$l$.
(Note that $\OT(\lambda,l)$ is empty unless $l-|\lambda|$ is nonnegative and even, where~$|\lambda|$ is the number of
squares in $\lambda$. Note also that if $l=|\lambda|$, then $\OT(\lambda,|\lambda|)$ is in simple bijection with the set
of standard Young tableaux of shape $\lambda$. In particular, the $(k+1)$th diagram of any $\eta\in\OT(\lambda,|\lambda|)$ is
obtained from the $k$th diagram of~$\eta$ by the addition of a square, for each $k=1,\ldots,|\lambda|$,
so if this square is in row~$i_k$ and column~$j_k$, then~$\eta$ can be
associated with the standard Young tableau of shape $\lambda$ in which the square in row~$i_k$ and column~$j_k$ is filled
with the integer~$k$.)
Now define an order $\prec$ on the integers by $\ldots\prec-2\prec2\prec-1\prec1\prec0$, and
let an ascent of an oscillating tableau~$\eta$ be an integer~$k$ for which $j_k-i_k\prec j_{k+1}-i_{k+1}$,
where~$i_k$ and~$j_k$ are the row and column of the square by which the~$k$th and $(k+1)$th diagrams of $\eta$ differ.
(Note that it is shown in~\cite{Behr-osc} that certain other orders on the integers can be used here instead of~$\prec$.)
Let $\mathrm{asc}(\eta)$ denote the number of ascents of~$\eta$.

For a strict partition $\kappa=(\kappa_1,\ldots,\kappa_r)$ of $p$ (i.e., a strictly-decreasing sequence of
positive integers whose sum is $p$), the double diagram~$\Delta(\kappa)$ is defined to be the Young diagram which is comprised
of~$r$ squares on the main diagonal, together with~$\kappa_i$ squares to the right of the main diagonal in row~$i$,
and~$\kappa_i-1$ squares below the main diagonal in column~$i$, for each $i=1,\ldots,r$.
(Using Frobenius notation,
$\Delta(\kappa)=(\kappa_1,\ldots,\kappa_r\mid \kappa_1-1,\ldots,\kappa_r-1)$.)
By using a bijection between ASMs and certain sets of osculating lattice paths,
a bijection between DPPs and certain sets of nonintersecting lattice paths,
and a result of Behrend~\cite[Cor.~14]{Behr-osc}, which itself follows from a bijection~\cite[Thm.~13]{Behr-osc} between
certain sets of paths and certain generalized oscillating tableaux, it can be shown that
\begin{align}\label{ASMosc}\hspace{4mm}|\{A\in\ASM(n)\mid\nu(A)=p\}|&=\sum_{\eta\in\OT(\emptyset,2p)}{n+\mathrm{asc}(\eta)\choose2p},\\
\label{DPPosc}|\{D\in\DPP(n)\mid\nu(D)=p\}|&=\sum_{\kappa\vDash p}\,\sum_{\eta\in\OT(\Delta(\kappa),2p)}{n+\mathrm{asc}(\eta)\choose2p},\end{align}
where $\kappa\vDash p$ denotes all strict partitions $\kappa$ of $p$.
(Note that~$|\Delta(\kappa)|=2p$, so $\OT(\Delta(\kappa),2p)$ is
in bijection with the set of standard Young tableaux of shape $\Delta(\kappa)$.)

As an example, consider $p=2$.  Then $\OT(\emptyset,4)$ consists of
\setlength{\unitlength}{3mm}
$(\emptyset,\begin{picture}(1.1,1)\multiput(0,0)(0,1){2}{\line(1,0){1}}\multiput(0,0)(1,0){2}{\line(0,1){1}}\end{picture},
\emptyset,\begin{picture}(1.1,1)\multiput(0,0)(0,1){2}{\line(1,0){1}}\multiput(0,0)(1,0){2}{\line(0,1){1}}\end{picture},\emptyset)$,
$(\emptyset,\begin{picture}(1.1,1)\multiput(0,0)(0,1){2}{\line(1,0){1}}\multiput(0,0)(1,0){2}{\line(0,1){1}}\end{picture},
\begin{picture}(2.1,1)\multiput(0,0)(0,1){2}{\line(1,0){2}}\multiput(0,0)(1,0){3}{\line(0,1){1}}\end{picture},
\begin{picture}(1.1,1)\multiput(0,0)(0,1){2}{\line(1,0){1}}\multiput(0,0)(1,0){2}{\line(0,1){1}}\end{picture},\emptyset)$ and
$\raisebox{-1mm}{$\Bigl($}\emptyset,\begin{picture}(1.1,1)\multiput(0,0)(0,1){2}{\line(1,0){1}}\multiput(0,0)(1,0){2}{\line(0,1){1}}\end{picture},
\raisebox{-3mm}{\begin{picture}(1.1,2)\multiput(0,0)(0,1){3}{\line(1,0){1}}\multiput(0,0)(1,0){2}{\line(0,1){2}}\end{picture}},
\begin{picture}(1.1,1)\multiput(0,0)(0,1){2}{\line(1,0){1}}\multiput(0,0)(1,0){2}{\line(0,1){1}}\end{picture},\emptyset\raisebox{-1mm}{$\Bigr)$}$,
which have~$0$,~$1$ and~$1$ ascent respectively. (For example, for the second of these oscillating tableaux,
the row $i_k$ and column $j_k$ of the square by which the~$k$th and $(k+1)$th diagrams differ are
$(i_1,j_1)=(i_4,j_4)=(1,1)$ and $(i_2,j_2)=(i_3,j_3)=(1,2)$, which gives the single ascent $3$,
since $j_3-i_3=1\prec0=j_4-i_4$.)  The only strict partition of $2$ is $(2)$, $\Delta(2)$ is
\raisebox{-3mm}{\begin{picture}(3.1,2)\multiput(0,1)(0,1){2}{\line(1,0){3}}\multiput(0,0)(1,0){2}{\line(0,1){2}}\put(0,0){\line(1,0){1}}
\multiput(2,1)(1,0){2}{\line(0,1){1}}\end{picture}}, and $\cup_{\kappa\vDash 2}\OT(\Delta(\kappa),4)=\OT(\Delta(2),4)$
consists of
$\raisebox{-1mm}{$\Bigl($}\emptyset,\begin{picture}(1.1,1)\multiput(0,0)(0,1){2}{\line(1,0){1}}\multiput(0,0)(1,0){2}{\line(0,1){1}}\end{picture},
\begin{picture}(2.1,1)\multiput(0,0)(0,1){2}{\line(1,0){2}}\multiput(0,0)(1,0){3}{\line(0,1){1}}\end{picture},
\raisebox{-3mm}{\begin{picture}(2.1,2)\multiput(0,1)(0,1){2}{\line(1,0){2}}\multiput(0,0)(1,0){2}{\line(0,1){2}}\put(0,0){\line(1,0){1}}
\put(2,1){\line(0,1){1}}\end{picture}},
\raisebox{-3mm}{\begin{picture}(3.1,2)\multiput(0,1)(0,1){2}{\line(1,0){3}}\multiput(0,0)(1,0){2}{\line(0,1){2}}\put(0,0){\line(1,0){1}}
\multiput(2,1)(1,0){2}{\line(0,1){1}}\end{picture}}\raisebox{-1mm}{$\Bigr)$}$,
$\raisebox{-1mm}{$\Bigl($}\emptyset,\begin{picture}(1.1,1)\multiput(0,0)(0,1){2}{\line(1,0){1}}\multiput(0,0)(1,0){2}{\line(0,1){1}}\end{picture},
\begin{picture}(2.1,1)\multiput(0,0)(0,1){2}{\line(1,0){2}}\multiput(0,0)(1,0){3}{\line(0,1){1}}\end{picture},
\begin{picture}(3.1,1)\multiput(0,0)(0,1){2}{\line(1,0){3}}\multiput(0,0)(1,0){4}{\line(0,1){1}}\end{picture},
\raisebox{-3mm}{\begin{picture}(3.1,2)\multiput(0,1)(0,1){2}{\line(1,0){3}}\multiput(0,0)(1,0){2}{\line(0,1){2}}\put(0,0){\line(1,0){1}}
\multiput(2,1)(1,0){2}{\line(0,1){1}}\end{picture}}\raisebox{-1mm}{$\Bigr)$}$ and
$\raisebox{-1mm}{$\Bigl($}\emptyset,\begin{picture}(1.1,1)\multiput(0,0)(0,1){2}{\line(1,0){1}}\multiput(0,0)(1,0){2}{\line(0,1){1}}\end{picture},
\raisebox{-3mm}{\begin{picture}(1.1,2)\multiput(0,0)(0,1){3}{\line(1,0){1}}\multiput(0,0)(1,0){2}{\line(0,1){2}}\end{picture}},
\raisebox{-3mm}{\begin{picture}(2.1,2)\multiput(0,1)(0,1){2}{\line(1,0){2}}\multiput(0,0)(1,0){2}{\line(0,1){2}}\put(0,0){\line(1,0){1}}
\put(2,1){\line(0,1){1}}\end{picture}},
\raisebox{-3mm}{\begin{picture}(3.1,2)\multiput(0,1)(0,1){2}{\line(1,0){3}}\multiput(0,0)(1,0){2}{\line(0,1){2}}\put(0,0){\line(1,0){1}}
\multiput(2,1)(1,0){2}{\line(0,1){1}}\end{picture}}\raisebox{-1mm}{$\Bigr)$}$, which have~$0$,~$1$ and~$1$ ascent respectively.
Therefore~\eqref{ASMosc} and \eqref{DPPosc} give $|\{A\in\ASM(n)\mid\nu(A)=2\}|=|\{D\in\ASM(n)\mid\nu(D)=2\}|={n\choose 4}+2{n+1\choose 4}$.
Some details of the further example of~\eqref{ASMosc} for $p=3$ can be found in~\cite[Sec.~17]{Behr-osc}.

It follows from~\eqref{ASMosc}--\eqref{DPPosc} that, for any $p$,
\begin{multline*}\hspace{4mm}|\{A\in\ASM(n)\mid\nu(A)=p\}|=|\{D\in\DPP(n)\mid\nu(D)=p\}|\text{ for all }n\text{ if and only if }\\
\ts|\{\eta\in\OT(\emptyset,2p)\mid\mathrm{asc}(\eta)=s\}|=|\{\eta\in\bigcup_{\kappa\vDash p}\OT(\Delta(\kappa),2p)\mid\mathrm{asc}(\eta)=s\}|
\text{ for all }s.\end{multline*}
Indeed, Theorem~\ref{MRRC1} of this paper implies that the first of these equalities between set sizes holds
for any~$n$ and~$p$, so that a corollary of Theorem~\ref{MRRC1} is
\begin{multline}\label{osceq}\hspace{4mm}|\{\eta\in\OT(\emptyset,2p)\mid\mathrm{asc}(\eta)=s\}|=\\
\ts|\{\eta\in\bigcup_{\kappa\vDash p}\OT(\Delta(\kappa),2p)\mid\mathrm{asc}(\eta)=s\}|
\text{ for any }p\text{ and }s.\end{multline}
Unfortunately, we do not have a direct or bijective proof of~\eqref{osceq}.
(Note that if a bijective proof were known then, due to the bijective nature of the
derivation of~\eqref{ASMosc} and~\eqref{DPPosc}, this would provide a bijection between
$\{A\in\ASM(n)\mid\nu(A)=p\}$ and $\{D\in\DPP(n)\mid\nu(D)=p\}$ for any~$n$ and~$p$.)
Nevertheless,~\eqref{osceq} can be proved bijectively for the special case
in which $s$ is summed over, i.e., an explicit bijection is known between
$\OT(\emptyset,2p)$ and $\cup_{\kappa\vDash p}\OT(\Delta(\kappa),2p)$ for any $p$
(but oscillating tableaux which correspond under this bijection do not always have
equal numbers of ascents).
This bijection is a composition of a bijection of Stanley (see
Roby~\cite[Sec.~4.2]{Roby-thesis}, Stanley~\cite[Sec.~9]{Stan-incrdecr},
Sundaram~\cite[Sec.~8]{Sunda-thesis} or Sundaram~\cite[Sec.~2]{Sunda}) between
$\OT(\emptyset,2p)$ and the set of fixed-point-free involutions on $\{1,\ldots,2p\}$
(or, equivalently, matchings of $\{1,\ldots,2p\}$,
or perfect matchings of the complete graph $K_{2p}$),
and a bijection, which follows from results of Burge~\cite[Sec.~4]{Burge},
between the set of such involutions and $\cup_{\kappa\vDash p}\OT(\Delta(\kappa),2p)$.
(It also follows that these sets all have size $(2p-1)!!$.)  We have also considered the further implications on~\eqref{osceq} provided by the
presence of the additional ASM and DPP statistics $\mu$ and $\rho$ in Theorem~\ref{MRRC1},
but this has so far not yielded a direct proof of~\eqref{osceq}.
\end{list}

\section{The unrefined enumeration ($z=1$)}\label{unref}
In this section, Theorem~\ref{MRRC1} will be proved for the case $z=1$ in~\eqref{MRRC2}.
This will be done by expressing the ASM and DPP generating functions~\eqref{ZASM} and~\eqref{ZDPP} for $z=1$ as
determinants of $n\times n$ matrices, and showing that these determinants are equal.

Let $\omega$ be a solution of
\begin{equation}\label{nuxy}y\omega^2+(1-x-y)\omega+x=0,\end{equation}
and define $n\times n$ matrices
\begin{align}\label{MASM}M_\ASM(n,x,y,1)_{ij}&=(1-\omega)\,\delta_{ij}+\omega\sum_{k=0}^{\min(i,j)}{i\choose k}{j\choose k}x^k\,y^{i-k},\\[1mm]
\label{MDPP}M_\DPP(n,x,y,1)_{ij}&=-\delta_{i,j+1}+\ds\sum_{k=0}^{\min(i,j+1)}{i-1\choose i-k}{j+1\choose k}x^k\,y^{i-k},\end{align}
with $0\le i,j\le n-1$. Note that $\omega$ depends on~$x$ and~$y$, even though this is not explicitly indicated by the notation.
Definitions of matrices $M_\ASM(n,x,y,z)$ and $M_\DPP(n,x,y,z)$ which reduce to~\eqref{MASM} and~\eqref{MDPP} for $z=1$ will be given
in Section~\ref{refenum}.

Note also that for $x=y=z=1$, which corresponds to the case of straight enumeration,
$M_\ASM(n,1,1,1)_{ij}=e^{\pm i\pi/3}\,\delta_{ij}+e^{\mp i\pi/3}{i+j\choose i}$ and
$M_\DPP(n,1,1,1)_{ij}=-\delta_{i,j+1}+{i+j\choose i}$.

Formulae for the $z=1$ ASM and DPP generating functions as determinants, and the equality of these determinants,
can now be stated as follows.
\begin{proposition}\label{ASMdet}Let $Z_\ASM(n,x,y,1)$ and $M_\ASM(n,x,y,1)$ be given by
\eqref{ZASM} and \eqref{MASM}.  Then
\[Z_\ASM(n,x,y,1)=\det M_\ASM(n,x,y,1).\]
\end{proposition}
\begin{proposition}\label{DPPdet} Let $Z_\DPP(n,x,y,1)$ and $M_\DPP(n,x,y,1)$ be given by
\eqref{ZDPP} and \eqref{MDPP}.  Then
\[Z_\DPP(n,x,y,1)=\det M_\DPP(n,x,y,1).\]
\end{proposition}
\begin{proposition}\label{ASMdetDPPdet}
Let $M_\ASM(n,x,y,1)$ and $M_\DPP(n,x,y,1)$ be given by~\eqref{MASM} and~\eqref{MDPP}. Then
\[\det M_\ASM(n,x,y,1)=\det M_\DPP(n,x,y,1).\]
\end{proposition}
The validity of Theorem~\ref{MRRC1} for $z=1$ in~\eqref{MRRC2}
follows immediately from Propositions~\ref{ASMdet}--\ref{ASMdetDPPdet}.
The proofs of these propositions will be given in the ensuing three subsections.

For further information on determinants of matrices closely related to those given here,
including contexts (usually that of plane partitions or rhombus tilings) in which they arise
and formulae for their evaluation in some cases, see for example
Andrews~\cite{And-weakMac,And-MacDPP}, Bressoud~\cite{Bressoud},
Ciucu, Eisenk{\"o}lbl, Krattenthaler and Zare~\cite[Thms.~10--13]{CEKZ},
Ciucu and Krattenthaler~\cite{CK-PP},
Colomo and Pronko~\cite[Eqs.~(23)--(24)]{CP-6v0},
Colomo and Pronko~\cite[Eqs.~(4.3)--(4.7)]{CP-6v},
Gessel and Xin~\cite[Sec.~5]{GX},
Krattenthaler~\cite[e.g., Thms.~25--37]{Kratt-det1},
Krattenthaler~\cite[Sec.~5.5]{Kratt-det2},
Mills, Robbins and Rumsey~\cite{MRR-Mac},
Mills, Robbins and Rumsey~\cite[Secs.~3--4]{MRR-symPP},
Mitra and Nienhuis~\cite[Sec.~5]{MN},
and Robbins~\cite[Sec.~2]{Rob-ASM}.

\subsection{\protect Proof of Proposition \ref{ASMdet}}
\label{ASMdetderiv}
In this subsection, the ASM determinant formula of Proposition \ref{ASMdet} will be obtained
using the Izergin--Korepin formula for the partition function of the
six-vertex model with domain-wall boundary conditions (DWBC), together
with a bijection between ASMs and configurations of this model.  Such a determinant formula
has also been obtained, using the Izergin--Korepin formula together with certain other results,
by Colomo and Pronko in~\cite[Eqs.~(23)--(24)]{CP-6v0} and~\cite[Eqs.~(4.3)--(4.7)]{CP-6v}.

The bijection between ASMs and configurations of the six-vertex model with DWBC will be outlined first.
This correspondence was first observed (using slightly different terminology) by Robbins and Rumsey~\cite[pp.~179--180]{RR-ASM}
and by Elkies, Kuperberg, Larsen and Propp~\cite[Sec.~7]{EKLP2}, and first applied to
the enumeration of ASMs by Kuperberg~\cite{Kup-ASM}.  Note that there is
a closely-related bijection between ASMs and certain sets of osculating lattice paths.
See, for example, Behrend~\cite[Secs.~2--4]{Behr-osc} for further information and references.
The statistical mechanical six-vertex or square ice model was itself first studied and solved by Lieb and Sutherland.
See, for example, Baxter~\cite[Chaps.~8~\&~9]{Baxter} for further information and references.
DWBC for this model were first introduced and studied by Korepin~\cite{Kor}.

Consider the $n\times n$ grid with vertices $\{(i,j)\mid i,j=0,\ldots,n+1\}\setminus\{(0,0),(0,n+1),(n+1,0),(n+1,n+1)\}$,
where $(i,j)$ is in the $i$th row from the top and $j$th column from the left, and for which there are
horizontal edges between $(i,j)$ and $(i,j\pm1)$, and vertical edges between $(i,j)$ and $(i\pm1,j)$,
for each $i,j=1,\ldots,n$.

\begin{figure}[h]
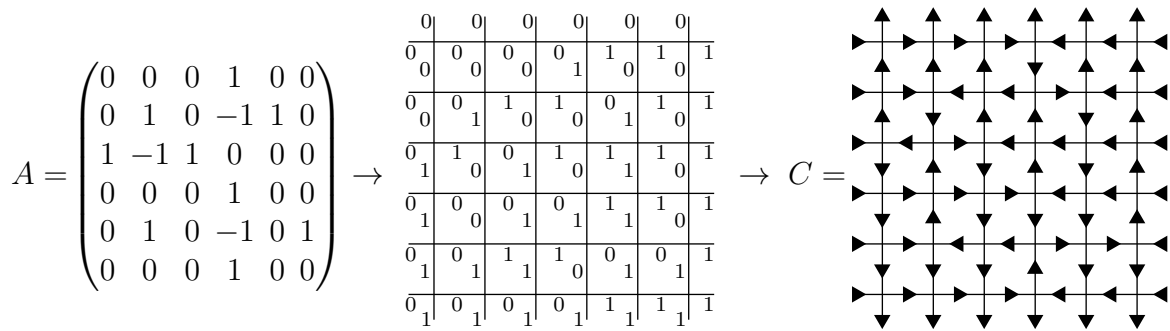
\centering
\raisebox{18mm}{$A=\!\left(\ba{c@{\,\;}c@{\,\;}c@{\,\;}c@{\,\;}c@{\,\;}c}0&0&0&1&0&0\\[0.2mm]0&1&0&-1&1&0\\[0.2mm]1&-1&1&0&0&0\\[0.2mm]0&0&0&1&0&0
\\[0.2mm]0&1&0&-1&0&1\\[0.2mm]0&0&0&1&0&0\ea
\right)\rightarrow\;$}
\psset{unit=6.7mm}
\pspicture(0.4,0.5)(6.8,6.5)
\multips(0.5,1)(0,1){6}{\psline[linewidth=0.5pt](0,0)(6,0)}\multips(1,0.5)(1,0){6}{\psline[linewidth=0.5pt](0,0)(0,6)}
\multirput[r](0.96,6.4)(1,0){6}{$\ss0$}\multirput[t](0.55,0.94)(0,1){6}{$\ss0$}\multirput[r](0.96,0.5)(1,0){6}{$\ss1$}
\multirput[t](6.4,0.94)(0,1){6}{$\ss1$}
\rput[t](1.46,0.94){$\ss0$}\rput[t](2.46,0.94){$\ss0$}\rput[t](3.46,0.94){$\ss0$}\rput[t](4.46,0.94){$\ss1$}\rput[t](5.46,0.94){$\ss1$}
\rput[t](1.46,1.94){$\ss0$}\rput[t](2.46,1.94){$\ss1$}\rput[t](3.46,1.94){$\ss1$}\rput[t](4.46,1.94){$\ss0$}\rput[t](5.46,1.94){$\ss0$}
\rput[t](1.46,2.94){$\ss0$}\rput[t](2.46,2.94){$\ss0$}\rput[t](3.46,2.94){$\ss0$}\rput[t](4.46,2.94){$\ss1$}\rput[t](5.46,2.94){$\ss1$}
\rput[t](1.46,3.94){$\ss1$}\rput[t](2.46,3.94){$\ss0$}\rput[t](3.46,3.94){$\ss1$}\rput[t](4.46,3.94){$\ss1$}\rput[t](5.46,3.94){$\ss1$}
\rput[t](1.46,4.94){$\ss0$}\rput[t](2.46,4.94){$\ss1$}\rput[t](3.46,4.94){$\ss1$}\rput[t](4.46,4.94){$\ss0$}\rput[t](5.46,4.94){$\ss1$}
\rput[t](1.46,5.94){$\ss0$}\rput[t](2.46,5.94){$\ss0$}\rput[t](3.46,5.94){$\ss0$}\rput[t](4.46,5.94){$\ss1$}\rput[t](5.46,5.94){$\ss1$}
\rput[r](0.94,1.46){$\ss1$}\rput[r](1.94,1.46){$\ss1$}\rput[r](2.94,1.46){$\ss1$}\rput[r](3.94,1.46){$\ss0$}\rput[r](4.94,1.46){$\ss1$}
\rput[r](5.94,1.46){$\ss1$}
\rput[r](0.94,2.46){$\ss1$}\rput[r](1.94,2.46){$\ss0$}\rput[r](2.94,2.46){$\ss1$}\rput[r](3.94,2.46){$\ss1$}\rput[r](4.94,2.46){$\ss1$}
\rput[r](5.94,2.46){$\ss0$}
\rput[r](0.94,3.46){$\ss1$}\rput[r](1.94,3.46){$\ss0$}\rput[r](2.94,3.46){$\ss1$}\rput[r](3.94,3.46){$\ss0$}\rput[r](4.94,3.46){$\ss1$}
\rput[r](5.94,3.46){$\ss0$}
\rput[r](0.94,4.46){$\ss0$}\rput[r](1.94,4.46){$\ss1$}\rput[r](2.94,4.46){$\ss0$}\rput[r](3.94,4.46){$\ss0$}\rput[r](4.94,4.46){$\ss1$}
\rput[r](5.94,4.46){$\ss0$}
\rput[r](0.94,5.46){$\ss0$}\rput[r](1.94,5.46){$\ss0$}\rput[r](2.94,5.46){$\ss0$}\rput[r](3.94,5.46){$\ss1$}\rput[r](4.94,5.46){$\ss0$}
\rput[r](5.94,5.46){$\ss0$}
\endpspicture
\raisebox{18mm}{$\;\rightarrow\;C=$}
\pspicture(0.5,0.5)(6.5,6.5)
\multips(0.5,1)(0,1){6}{\psline[linewidth=0.5pt](0,0)(6,0)}\multips(1,0.5)(1,0){6}{\psline[linewidth=0.5pt](0,0)(0,6)}
\psdots[dotstyle=triangle*,dotscale=1.5,dotangle=0]
(1,6.5)(2,6.5)(3,6.5)(4,6.5)(5,6.5)(6,6.5)
(1,5.5)(2,5.5)(3,5.5)(5,5.5)(6,5.5)
(1,4.5)(3,4.5)(4,4.5)(6,4.5)
(2,3.5)(4,3.5)(6,3.5)
(2,2.5)(6,2.5)
(4,1.5)
\psdots[dotstyle=triangle*,dotscale=1.5,dotangle=180]
(1,0.5)(2,0.5)(3,0.5)(4,0.5)(5,0.5)(6,0.5)
(1,1.5)(2,1.5)(3,1.5)(5,1.5)(6,1.5)
(1,2.5)(3,2.5)(4,2.5)(5,2.5)
(1,3.5)(3,3.5)(5,3.5)
(2,4.5)(5,4.5)
(4,5.5)
\psdots[dotstyle=triangle*,dotscale=1.5,dotangle=90]
(6.5,1)(6.5,2)(6.5,3)(6.5,4)(6.5,5)(6.5,6)
(5.5,1)(5.5,3)(5.5,4)(5.5,5)(5.5,6)
(4.5,1)(4.5,3)(4.5,4)(4.5,6)
(3.5,2)(3.5,4)(3.5,5)
(2.5,2)(2.5,5)
(1.5,4)
\psdots[dotstyle=triangle*,dotscale=1.5,dotangle=-90]
(0.5,1)(0.5,2)(0.5,3)(0.5,4)(0.5,5)(0.5,6)
(1.5,1)(1.5,2)(1.5,3)(1.5,5)(1.5,6)
(2.5,1)(2.5,3)(2.5,4)(2.5,6)
(3.5,1)(3.5,3)(3.5,6)
(4.5,2)(4.5,5)
(5.5,2)
\endpspicture
\caption{The mapping from an ASM to a six-vertex model configuration.}
\label{ASMPEx}\end{figure}

Now let $\SVDWBC(n)$ be the set of all configurations of the six-vertex model on the $n\times n$ grid with DWBC, i.e.,
decorations of the grid's edges with arrows such that:
\begin{itemize}
\item The arrows on the external edges are fixed, with the horizontal ones all incoming  and the
vertical ones all outgoing.
\item At each internal vertex, there are two incoming and two outgoing arrows.
\end{itemize}
The latter condition is the ``six-vertex'' condition, since it allows for
only six possible arrow configurations around an internal vertex, denoted by
the symbols $a_1,a_2,b_1,b_2,c_1,c_2$:
\psset{unit=14mm}\begin{center}\pspicture(0.5,-0.4)(11.5,1)
\multips(1,0.5)(2,0){6}{\vertex}
\rput(1,-0.3){$a_1$}
\psdots[dotstyle=triangle*,dotscale=1.5](1,0.75)
\psdots[dotstyle=triangle*,dotscale=1.5](1,0.25)
\psdots[dotstyle=triangle*,dotscale=1.5,dotangle=-90](0.75,0.5)
\psdots[dotstyle=triangle*,dotscale=1.5,dotangle=-90](1.25,0.5)
\rput(3,-0.3){$a_2$}
\psdots[dotstyle=triangle*,dotscale=1.5,dotangle=180](3,0.75)
\psdots[dotstyle=triangle*,dotscale=1.5,dotangle=180](3,0.25)
\psdots[dotstyle=triangle*,dotscale=1.5,dotangle=90](2.75,0.5)
\psdots[dotstyle=triangle*,dotscale=1.5,dotangle=90](3.25,0.5)
\rput(5,-0.3){$b_1$}
\psdots[dotstyle=triangle*,dotscale=1.5](5,0.75)
\psdots[dotstyle=triangle*,dotscale=1.5](5,0.25)
\psdots[dotstyle=triangle*,dotscale=1.5,dotangle=90](4.75,0.5)
\psdots[dotstyle=triangle*,dotscale=1.5,dotangle=90](5.25,0.5)
\rput(7,-0.3){$b_2$}
\psdots[dotstyle=triangle*,dotscale=1.5,dotangle=180](7,0.75)
\psdots[dotstyle=triangle*,dotscale=1.5,dotangle=180](7,0.25)
\psdots[dotstyle=triangle*,dotscale=1.5,dotangle=-90](6.75,0.5)
\psdots[dotstyle=triangle*,dotscale=1.5,dotangle=-90](7.25,0.5)
\rput(9,-0.3){$c_1$}
\psdots[dotstyle=triangle*,dotscale=1.5](9,0.75)
\psdots[dotstyle=triangle*,dotscale=1.5,dotangle=180](9,0.25)
\psdots[dotstyle=triangle*,dotscale=1.5,dotangle=-90](8.75,0.5)
\psdots[dotstyle=triangle*,dotscale=1.5,dotangle=90](9.25,0.5)
\rput(11,-0.3){$c_2$}
\psdots[dotstyle=triangle*,dotscale=1.5,dotangle=180](11,0.75)
\psdots[dotstyle=triangle*,dotscale=1.5](11,0.25)
\psdots[dotstyle=triangle*,dotscale=1.5,dotangle=90](10.75,0.5)
\psdots[dotstyle=triangle*,dotscale=1.5,dotangle=-90](11.25,0.5)
\endpspicture\end{center}
The mappings between $\ASM(n)$ and $\SVDWBC(n)$,
which can be shown straightforwardly to be well-defined and bijective, are as follows.
\begin{itemize}
\item To map $C\in\SVDWBC(n)$ to $A\in\ASM(n)$, place a~$0$ at vertex configurations in $C$ of types~$a_1$,
$a_2$, $b_1$ or $b_2$, a~$1$ at configurations of type $c_1$ and a~$-1$ at configurations of type~$c_2$.
This array of $0$'s, $1$'s and $-1$'s is then the ASM $A$.
\item To map $A\in\ASM(n)$ to $C\in\SVDWBC(n)$, first associate the partial row sum $\sum_{j'=1}^{j}A_{ij'}$
with the horizontal edge between $(i,j)$ and $(i,j+1)$,
for each $i=1,\ldots,n$, $j=0,\ldots,n$, and associate
the partial column sum $\sum_{i'=1}^{i}A_{i'j}$ with the vertical edge between $(i,j)$ and $(i+1,j)$,
for each $i=0,\ldots,n$, $j=1,\ldots,n$.
The defining properties of ASMs ensure that each of these partial sums is~0 or~1.  Now
obtain $C$ by placing a right/up (respectively left/down) arrow on each edge associated with a $0$ (respectively $1$). This mapping applied to
the example of~\eqref{ASMEx} is shown in Figure~\ref{ASMPEx}, where the intermediate grid shows the values of the partial sums
associated with each edge.
\end{itemize}

As a further example, mapping $\ASM(3)$, as given in~\eqref{ASM3}, to $\SVDWBC(3)$ gives
\psset{unit=4.4mm}
\begin{equation}\label{6VDBWC3}\SVDWBC(3)=\left\{\raisebox{-5mm}{
\pspicture(0.5,0.5)(4.4,3.5)\multips(0.5,1)(0,1){3}{\psline[linewidth=0.5pt](0,0)(3,0)}\multips(1,0.5)(1,0){3}{\psline[linewidth=0.5pt](0,0)(0,3)}
\psdots[dotstyle=triangle*,dotscale=1,dotangle=0](1,3.5)(2,2.5)(2,3.5)(3,1.5)(3,2.5)(3,3.5)
\psdots[dotstyle=triangle*,dotscale=1,dotangle=180](1,0.5)(1,1.5)(1,2.5)(2,0.5)(2,1.5)(3,0.5)
\psdots[dotstyle=triangle*,dotscale=1,dotangle=90](1.5,3)(2.5,2)(2.5,3)(3.5,1)(3.5,2)(3.5,3)
\psdots[dotstyle=triangle*,dotscale=1,dotangle=-90](0.5,1)(0.5,2)(0.5,3)(1.5,1)(1.5,2)(2.5,1)\rput(3.9,1.2){,}\endpspicture
\pspicture(0.5,0.5)(4.4,3.5)\multips(0.5,1)(0,1){3}{\psline[linewidth=0.5pt](0,0)(3,0)}\multips(1,0.5)(1,0){3}{\psline[linewidth=0.5pt](0,0)(0,3)}
\psdots[dotstyle=triangle*,dotscale=1,dotangle=0](1,1.5)(1,2.5)(1,3.5)(2,2.5)(2,3.5)(3,3.5)
\psdots[dotstyle=triangle*,dotscale=1,dotangle=180](1,0.5)(2,0.5)(2,1.5)(3,0.5)(3,1.5)(3,2.5)
\psdots[dotstyle=triangle*,dotscale=1,dotangle=90](1.5,1)(2.5,1)(2.5,2)(3.5,1)(3.5,2)(3.5,3)
\psdots[dotstyle=triangle*,dotscale=1,dotangle=-90](0.5,1)(0.5,2)(0.5,3)(1.5,2)(1.5,3)(2.5,3)\rput(3.9,1.2){,}\endpspicture
\pspicture(0.5,0.5)(4.4,3.5)\multips(0.5,1)(0,1){3}{\psline[linewidth=0.5pt](0,0)(3,0)}\multips(1,0.5)(1,0){3}{\psline[linewidth=0.5pt](0,0)(0,3)}
\psdots[dotstyle=triangle*,dotscale=1,dotangle=0](1,3.5)(2,1.5)(2,2.5)(2,3.5)(3,2.5)(3,3.5)
\psdots[dotstyle=triangle*,dotscale=1,dotangle=180](1,0.5)(1,1.5)(1,2.5)(2,0.5)(3,0.5)(3,1.5)
\psdots[dotstyle=triangle*,dotscale=1,dotangle=90](1.5,3)(2.5,1)(2.5,3)(3.5,1)(3.5,2)(3.5,3)
\psdots[dotstyle=triangle*,dotscale=1,dotangle=-90](0.5,1)(0.5,2)(0.5,3)(1.5,1)(1.5,2)(2.5,2)\rput(3.9,1.2){,}\endpspicture
\pspicture(0.5,0.5)(4.4,3.5)\multips(0.5,1)(0,1){3}{\psline[linewidth=0.5pt](0,0)(3,0)}\multips(1,0.5)(1,0){3}{\psline[linewidth=0.5pt](0,0)(0,3)}
\psdots[dotstyle=triangle*,dotscale=1,dotangle=0](1,2.5)(1,3.5)(2,1.5)(2,2.5)(2,3.5)(3,3.5)
\psdots[dotstyle=triangle*,dotscale=1,dotangle=180](1,0.5)(1,1.5)(2,0.5)(3,0.5)(3,1.5)(3,2.5)
\psdots[dotstyle=triangle*,dotscale=1,dotangle=90](1.5,2)(2.5,1)(2.5,2)(3.5,1)(3.5,2)(3.5,3)
\psdots[dotstyle=triangle*,dotscale=1,dotangle=-90](0.5,1)(0.5,2)(0.5,3)(1.5,1)(1.5,3)(2.5,3)\rput(3.9,1.2){,}\endpspicture
\pspicture(0.5,0.5)(4.4,3.5)\multips(0.5,1)(0,1){3}{\psline[linewidth=0.5pt](0,0)(3,0)}\multips(1,0.5)(1,0){3}{\psline[linewidth=0.5pt](0,0)(0,3)}
\psdots[dotstyle=triangle*,dotscale=1,dotangle=0](1,2.5)(1,3.5)(2,3.5)(3,1.5)(3,2.5)(3,3.5)
\psdots[dotstyle=triangle*,dotscale=1,dotangle=180](1,0.5)(1,1.5)(2,0.5)(2,1.5)(2,2.5)(3,0.5)
\psdots[dotstyle=triangle*,dotscale=1,dotangle=90](1.5,2)(2.5,2)(2.5,3)(3.5,1)(3.5,2)(3.5,3)
\psdots[dotstyle=triangle*,dotscale=1,dotangle=-90](0.5,1)(0.5,2)(0.5,3)(1.5,1)(1.5,3)(2.5,1)\rput(3.9,1.2){,}\endpspicture
\pspicture(0.5,0.5)(4.4,3.5)\multips(0.5,1)(0,1){3}{\psline[linewidth=0.5pt](0,0)(3,0)}\multips(1,0.5)(1,0){3}{\psline[linewidth=0.5pt](0,0)(0,3)}
\psdots[dotstyle=triangle*,dotscale=1,dotangle=0](1,1.5)(1,2.5)(1,3.5)(2,3.5)(3,2.5)(3,3.5)
\psdots[dotstyle=triangle*,dotscale=1,dotangle=180](1,0.5)(2,0.5)(2,1.5)(2,2.5)(3,0.5)(3,1.5)
\psdots[dotstyle=triangle*,dotscale=1,dotangle=90](1.5,1)(2.5,1)(2.5,3)(3.5,1)(3.5,2)(3.5,3)
\psdots[dotstyle=triangle*,dotscale=1,dotangle=-90](0.5,1)(0.5,2)(0.5,3)(1.5,2)(1.5,3)(2.5,2)\rput(3.9,1.2){,}\endpspicture
\pspicture(0.5,0.5)(3.7,3.5)\multips(0.5,1)(0,1){3}{\psline[linewidth=0.5pt](0,0)(3,0)}\multips(1,0.5)(1,0){3}{\psline[linewidth=0.5pt](0,0)(0,3)}
\psdots[dotstyle=triangle*,dotscale=1,dotangle=0](1,2.5)(1,3.5)(2,1.5)(2,3.5)(3,2.5)(3,3.5)
\psdots[dotstyle=triangle*,dotscale=1,dotangle=180](1,0.5)(1,1.5)(2,0.5)(2,2.5)(3,0.5)(3,1.5)
\psdots[dotstyle=triangle*,dotscale=1,dotangle=90](1.5,2)(2.5,1)(2.5,3)(3.5,1)(3.5,2)(3.5,3)
\psdots[dotstyle=triangle*,dotscale=1,dotangle=-90](0.5,1)(0.5,2)(0.5,3)(1.5,1)(1.5,3)(2.5,2)\endpspicture}
\right\},\end{equation}
where the elements are listed in the order corresponding to that used in~\eqref{ASM3}.

For $C\in\SVDWBC(n)$, denote the total number of type-$t$ vertex configurations in~$C$ as~$N_t(C)$,
for $t\in\{a_1,a_2,b_1,b_2,c_1,c_2\}$.
It can be seen that the arrow on the vertical edge immediately above $(i,j)$
and the arrow on the horizontal edge immediately to the right of $(i,j)$ both point towards $(i,j)$
if and only if a type-$a_2$ configuration occurs in $C$ at vertex $(i,j)$.
Now let $A\in\ASM(n)$ correspond to $C$. It then follows from the previous observation, and from the
mapping from $\ASM(n)$ to $\SVDWBC(n)$,
that $\sum_{i'=1}^{i-1}A_{i'j}\sum_{j'=1}^jA_{ij'}$ is 1
if a type-$a_2$ configuration occurs at $(i,j)$, or~$0$ if any of the other five configuration types occurs at $(i,j)$.
Therefore,  $N_{a_2}(C)=\sum_{i,j=1}^n\sum_{j'=1}^j\sum_{i'=1}^{i-1}A_{ij'}A_{i'j}
=\sum_{1\le i'<i\le n,\:1\le j'\le j\le n}A_{ij'}\,A_{i'j}=\nu(A)$, where~$\nu$ is the ASM
statistic given by~\eqref{nuA}.  It can be shown similarly that $N_{a_1}(C)=\nu(A)$ and
$N_{b_2}(C)=N_{b_1}(C)=\sum_{1\le i<i'\le n,\:1\le j\le j'\le n}A_{ij}\,A_{i'j'}=\sum_{1\le i\le i'\le n,\:1\le j<j'\le n}A_{ij}\,A_{i'j'}$.
It also follows from the bijection between~$\ASM(n)$ and~$\SVDWBC(n)$,
and the fact that each row and column of an ASM contains one more~1 than~$-1$, that
$N_{c_2}(C)=N_{c_1}(C)-n=\mu(A)$, where~$\mu$ is the ASM statistic given by~\eqref{muA}.

It is therefore natural to define, for any $C\in\SVDWBC(n)$,
\begin{align}\label{Na1a2}N_a(C)&=N_{a_1}(C)=N_{a_2}(C),\\
N_b(C)&=N_{b_1}(C)=N_{b_2}(C),\\
\label{Nc1c2}N_c(C)&=N_{c_1}(C)-n=N_{c_2}(C).\end{align}
Since $\sum_{t\in\{a_1,a_2,b_1,b_2,c_1,c_2\}}N_t(C)=n^2$, these numbers satisfy
\begin{equation}\label{paramsum}N_a(C)+N_b(C)+N_c(C)=n(n-1)/2.\end{equation}
Also, the previous observations of their relation with the ASM statistics~\eqref{nuA} and~\eqref{muA} can be summarized as follows.
\begin{lemma}\label{numuASM}Let $A\in \ASM(n)$ correspond to $C\in \SVDWBC(n)$.  Then
\begin{align*}\nu(A)&=N_a(C),\\
\mu(A)&=N_c(C).\end{align*}\end{lemma}

Now let $u$ and $v$ be parameters (often known as spectral parameters),
and associate a weight $\bar{a}(u,v)$, $\bar{b}(u,v)$ or $\bar{c}(u,v)$ to arrow configurations around a vertex
of types $a_1$ and~$a_2$,~$b_1$ and $b_2$ or $c_1$ and $c_2$ respectively, where
\begin{equation}\label{wuv}\ds\bar{a}(u,v)=uq-\frac{1}{vq},\qquad\bar{b}(u,v)=\frac{u}{q}-\frac{q}{v},\qquad
\bar{c}(u,v)=\Bigl(q^2-\frac{1}{q^2}\Bigr)
\sqrt{\frac{u}{v}},\end{equation}
and~$q$ is a further fixed parameter.
The six-vertex model
weights take various forms in the literature, with the form of~\eqref{wuv} seeming the most appropriate for the
purposes of this paper.  These weights satisfy $\bigl(\bar{a}(u,v)^2+\bar{b}(u,v)^2-\bar{c}(u,v)^2\bigr)\big/\bigl(
\bar{a}(u,v)\bar{b}(u,v)\bigr)=q^2+q^{-2}$, and the fact that
this is independent of $u$ and $v$ implies (see for example Baxter~\cite[pp.~187--190]{Baxter}) integrability, in the sense
that the Yang--Baxter equation is satisfied.
The form of the Yang--Baxter equation satisfied by the weights of~\eqref{wuv} is
\setlength{\unitlength}{10mm}
\begin{equation}\label{YBE}\raisebox{-18mm}{\begin{picture}(6,4)\put(0,1){\line(1,1){2}}\put(0,3){\line(1,-1){2}}
\put(3,0){\line(0,1){4}}\multiput(2,1)(0,2){2}{\line(1,0){2}}
\put(0.85,2){\pictss{r}{q w_1,\frac{q}{w_2}}}\put(2.95,0.9){\pictss{tr}{w_1,w_3}}\put(2.95,2.9){\pictss{tr}{w_2,w_3}}
\put(5,2){\pict{}{=}}\end{picture}
\begin{picture}(4,4)\put(1,0){\line(0,1){4}}\multiput(0,1)(0,2){2}{\line(1,0){2}}
\put(2.85,2){\pictss{r}{q w_1,\frac{q}{w_2}}}\put(0.95,0.9){\pictss{tr}{w_2,w_3}}\put(0.95,2.9){\pictss{tr}{w_1,w_3}}
\put(2,1){\line(1,1){2}}\put(2,3){\line(1,-1){2}}\end{picture}}\end{equation}
where this holds for any arrow configuration on the six external edges and any $w_1$, $w_2$ and~$w_3$, and there is
summation over all arrow configurations on the three internal edges of each side consistent with the six-vertex condition
at each vertex.

Now associate parameters~$u_1,\ldots,u_n$ with the rows, and parameters~$v_1,\ldots,v_n$ with the columns of the $n\times n$ grid.
For $C\in\SVDWBC(n)$, assign a weight $W_{ij}(C)$ which depends on~$u_i$ and~$v_j$ to
the configuration of arrows of $C$ around vertex $(i,j)$.  The partition function for this model is
\begin{equation}\label{ZSVDWBC}Z_\SVDWBC(u_1,\ldots,u_n;v_1,\ldots,v_n)=
\sum_{C\in\SVDWBC(n)}\;\prod_{i,j=1}^n\,W_{ij}(C).\end{equation}
It was shown by Izergin~\cite{Iz-6v}, using certain results of Korepin~\cite{Kor}, that
this partition function, for weights which satisfy a Yang--Baxter equation, can be expressed in terms of the determinant of an $n\times n$ matrix.
For weights given by~\eqref{wuv}, this Izergin--Korepin determinant formula is as follows.
\begin{theorem*}[Izergin]
\begin{multline}\label{Idet1}
Z_\SVDWBC(u_1,\ldots,u_n;v_1,\ldots,v_n)=\\
\frac{\prod_{i=1}^n\bar{c}(u_i,v_i)\:\prod_{i,j=1}^n\bar{a}(u_i,v_j)\bar{b}(u_i,v_j)}{\prod_{1\le i<j\le n}(u_i-u_j)(v_j^{-1}-v_i^{-1})}
\det_{1\le i,j\le n}\left(\frac{1}{\bar{a}(u_i,v_j)\bar{b}(u_i,v_j)}\right).\end{multline}
\end{theorem*}
This can be proved by showing that each side satisfies, and is uniquely determined by, four particular properties.
For example, one of these properties is that each side of~\eqref{Idet1} is symmetric in
$u_1,\ldots,u_n$ and (separately) in $v_1,\ldots,v_n$.  This is immediate for the RHS, and can be obtained for the LHS using the Yang--Baxter
equation~\eqref{YBE}.  For details of this proof of~\eqref{Idet1} see Izergin~\cite{Iz-6v}, Izergin, Coker and Korepin~\cite{ICK} or
Korepin and Zinn-Justin~\cite{artic12}.
For an alternative proof see Bogoliubov, Pronko and Zvonarev~\cite[Sec.~4]{BPZ-6v}.

Using \eqref{wuv}, \eqref{Idet1} can also be written as
\begin{multline}Z_\SVDWBC(u_1,\ldots,u_n;v_1,\ldots,v_n)\\
=\frac{\prod_{i=1}^nu_i^{1/2}v_i^{n+1/2}\:\prod_{i,j=1}^n\bar{a}(u_i,v_j)\bar{b}(u_i,v_j)}{\prod_{1\le i<j\le n}(u_i-u_j)(v_i-v_j)}
\det_{1\le i,j\le n}\left(\frac{1}{u_iv_j-q^2}-\frac{1}{u_iv_j-q^{-2}}\right).
\label{Idet2}
\end{multline}
Now consider the homogeneous case in which $u_1=\cdots=u_n=v_1=\cdots=v_n=r$, for a parameter $r$, and define $a$, $b$ and $c$ as
\begin{equation}\label{wr}\ds a=\bar{a}(r,r)=qr-(qr)^{-1},\quad b=\bar{b}(r,r)=q^{-1}r-qr^{-1},\quad
c=\bar{c}(r,r)=q^2-q^{-2}.\end{equation}
Note that $a$, $b$ and $c$ will remain dependent on $q$ and $r$, although this is not indicated explicitly by the notation.
The partition function~\eqref{ZSVDWBC} with these weights is, using~\eqref{Na1a2}--\eqref{Nc1c2}, given by
\begin{equation}Z_\SVDWBC(r,\ldots,r;r,\ldots,r)=\sum_{C\in\SVDWBC(n)}a^{2N_a(C)}\,b^{2N_b(C)}\,c^{2N_c(C)+n}.\end{equation}
The bijection between~$\ASM(n)$ and~$\SVDWBC(n)$, together with \eqref{ZASM}, \eqref{paramsum} and Lemma~\ref{numuASM}, now gives
\begin{equation}\label{ZCZASM}\ts Z_\SVDWBC(r,\ldots,r;r,\ldots,r)=
b^{n(n-1)}c^n\,Z_\ASM\bigl(n,\bigl(\frac{a}{b}\bigr)^2,\bigl(\frac{c}{b}\bigr)^2,1\bigr).\end{equation}
Therefore, $r$ and $q$ can be regarded as parameterizing $x$ and $y$ in~\eqref{ZASM}, where $x=\bigl(\frac{a}{b}\bigr)^2=
\bigl(\frac{qr-(qr)^{-1}}{q^{-1}r-qr^{-1}}\bigr)^2$
and $y=\bigl(\frac{c}{b}\bigr)^2=\bigl(\frac{q^2-q^{-2}}{q^{-1}r-qr^{-1}}\bigr)^2$.

An expression for $Z_\SVDWBC(r,\ldots,r;r,\ldots,r)$ cannot be obtained immediately from~\eqref{Idet2},
since the denominator of the prefactor and the determinant both vanish.  However, such an expression can be obtained
using the following general results.

Consider a power series $f(u,v)$ and parameters $u_0,\ldots,u_{n-1}$, $v_0,\ldots,v_{n-1}$.  Then
\begin{equation}\label{DDdet}
\frac{\ds\det_{0\le i,j\le n-1}\Bigl(f(u_i,v_j)\Bigr)}{\prod_{0\le i<j\le n-1}(u_j-u_i)(v_j-v_i)}
=\det_{0\le i,j\le n-1}\Bigl(f[u_0,\ldots,u_i][v_0,\ldots,v_j]\Bigr),\end{equation}
where $f[u_0,\ldots,u_i][v_0,\ldots,v_j]$ is the divided difference,
\begin{equation}\label{DD1}f[u_0,\ldots,u_i][v_0,\ldots,v_j]=
\sum_{k=0}^i\sum_{l=0}^j\frac{f(u_k,v_l)}{\prod_{\stackrel{k'=0}{\sss k'\ne k}}^i(u_k-u_{k'})
\prod_{\stackrel{l'=0}{\sss l'\ne l}}^j(v_l-v_{l'})}\,.\end{equation}
This result can be obtained by considering the matrix
\[D(u_0,\ldots,u_{n-1})_{ij}=\begin{cases}\prod_{\stackrel{k=0}{\sss k\ne j}}^i(u_j-u_k)^{-1},&i\ge j\\[2mm]0,&i<j\end{cases}\]
with $0\le i,j\le n-1$.
Then $\det D(u_0,\ldots,u_{n-1})=\prod_{0\le i<j\le n-1}(u_j-u_i)^{-1}$, and the matrix on the RHS of~\eqref{DDdet} is related to
the matrix on the LHS of~\eqref{DDdet} by multiplication on the left by $D(u_0,\ldots,u_{n-1})$ and on the right by
$D(v_0,\ldots,v_{n-1})^t$.

Let $[u^iv^j]f(u,v)$ denote the coefficient of~$u^iv^j$ in a power series~$f(u,v)$.
It can be shown that the divided difference~\eqref{DD1} can be written as
\begin{equation}f[u_0,\ldots,u_i][v_0,\ldots,v_j]=\sum_{k,l=0}^\infty\,d_{i+k,j+l}\,h_k(u_0,\ldots,u_i)\,h_l(v_0,\ldots,v_j),\end{equation}
where $d_{ij}=[u^iv^j]f(u,v)$, and  $h_k(u_0,\ldots,u_i)$ is the complete symmetric polynomial, defined by
\begin{equation}h_k(u_0,\ldots,u_i)=\begin{cases}\;1,&k=0\\[1mm]
\ds\sum_{0\le m_1\le\ldots\le m_k\le i}u_{m_1}\ldots u_{m_k},&k>0.\end{cases}\end{equation}
It now follows that the RHS of~\eqref{DDdet} is well-defined for
$u_0=\cdots=u_{n-1}$ and $v_0=\cdots=v_{n-1}$.  Furthermore, setting $r=u_0=\cdots=u_{n-1}$ and $s=v_0=\cdots=v_{n-1}$ gives
\begin{align}\label{DD}
f[\,\underbrace{r,\ldots,r}_{i+1}\,][\,\underbrace{s,\ldots,s}_{j+1}\,]&=
\sum_{k,l=0}^\infty\,d_{i+k,j+l}\,{i+k\choose k}{j+l\choose l}\,r^k\,s^l\\[2mm]
&=[u^iv^j]f(u+r,v+s).\nonumber
\end{align}
The homogeneous limit of~\eqref{Idet2} can now be obtained, using \eqref{DDdet} and \eqref{DD}, as
\begin{multline}\label{Idethomog}
Z_\SVDWBC(r,\ldots,r;r,\ldots,r)=\\
r^{n(n+1)}\:(ab)^{n^2}\!\ds\det_{0\le i,j\le n-1}\ts\left(
[u^iv^j]\left(\frac{1}{(u+r)(v+r)-q^2}-\frac{1}{(u+r)(v+r)-q^{-2}}\right)\right).\end{multline}
A homogeneous limit of an Izergin--Korepin determinant was first obtained by Izergin, Coker and Korepin in~\cite[Sec.~7]{ICK}.

Now define, for parameters $\alpha$ and $\beta$, the $n\times n$ lower triangular matrix
\begin{align}\label{L}L(\alpha,\beta)_{ij}&={i\choose j}\alpha^i\beta^j\\
\notag&=[u^iv^j]\,\frac{1}{1-\alpha u(1+\beta v)},\end{align}
with $0\le i,j\le n-1$.
Then the following properties are satisfied:
\begin{align}\label{LLT}(L(\alpha,\beta)\,L(\alpha,\beta)^t)_{ij}&=\ts\sum_{k=0}^{\min(i,j)}{i\choose k}{j\choose k}\alpha^{i+j}\,\beta^{2k},\\[2mm]
\label{LLTGF}(L(\alpha,\beta)\,L(\alpha,\beta)^t)_{ij}
&\ts=[u^iv^j]\,\frac{1}{1-\alpha(u+v)-\alpha^2(\beta^2-1)uv},\\[2mm]
\label{Ldet}\det L(\alpha,\beta)&=(\alpha\beta)^{n(n-1)/2},\\[2mm]
\label{LL}\ts L(\alpha_1,\beta_1)\,L(\frac{\alpha_2-\alpha_1}{\alpha_1\beta_1},\frac{\alpha_2\beta_2}{\alpha_2-\alpha_1})&=L(\alpha_2,\beta_2),
\text{ for any $\alpha_1$, $\alpha_2$, $\beta_1$, $\beta_2$}.
\end{align}
Properties \eqref{LLT} and \eqref{Ldet} follow immediately from~\eqref{L}, while properties~\eqref{LLTGF}
and~\eqref{LL} can be obtained easily using Cauchy-type residue integrals with appropriately-chosen contours
of integration, or using certain other standard methods for proving binomial coefficient identities.
For example, \eqref{LL} follows from
\begin{align}\label{LLderiv}
\bigl(L(\alpha_1,\beta_1)\,L(\alpha,\beta)\bigr)_{ij}
&=[u^iv^j]\oint\frac{dw}{2\pi iw}\,\frac{1}{1-\alpha_1u(1+\beta_1w)}\,\frac{1}{1-\frac{\alpha}{w}(1+\beta v)}\\
\notag&=[u^iv^j]\oint\frac{dw}{2\pi i}\,\frac{1}{1-\alpha_1u(1+\beta_1 w)}\,\frac{1}{w-{\alpha}(1+\beta v)}\\
\notag&=[u^iv^j]\frac{1}{1-\alpha_1 u(1+\beta_1 \alpha(1+\beta v))}\\
\notag&=L(\alpha_2,\beta_2)_{ij},
\end{align}
where $\alpha_2=\alpha_1(1+\beta_1\alpha)$ and $\alpha_2\beta_2=\alpha_1\beta_1\alpha\beta$.

Also define
\begin{equation}\label{nuqr}\omega_{\pm}=\frac{r^{\pm2}-q^{\mp2}}{q^{\pm2}-q^{\mp2}},\end{equation}
so that
\begin{equation}\label{nuabc}c^2\omega_\pm^2+(b^2-a^2-c^2)\omega_\pm+a^2=0.\end{equation}
We now have
\begin{alignat*}{2}
&\ts Z_\ASM\bigl(n,\bigl(\frac{a}{b}\bigr)^2,\bigl(\frac{c}{b}\bigr)^2,1\bigr)=
b^{-n(n-1)}\,c^{-n}\,Z_\SVDWBC(r,\ldots,r;r,\ldots,r)&\llap{[using \eqref{ZCZASM}]}\\
&\ts=\frac{r^{n(n+1)}\,a^{n^2}b^n}{c^n}\ds\det_{0\le i,j\le n-1}\ts\left(
[u^iv^j]\left(\frac{1}{(u+r)(v+r)-q^2}-\frac{1}{(u+r)(v+r)-q^{-2}}\right)\right)&\llap{\text{[using \eqref{Idethomog}]}}\\
&\ts=\frac{r^{n(n+1)}\,a^{n^2}b^n}{c^n}\,\det\Bigl(\frac{1}{r^2-q^2}\,L\bigl(\frac{r}{q^2-r^2},\frac{q}{r}\bigr)\,
L\bigl(\frac{r}{q^2-r^2},\frac{q}{r}\bigr)^t\,-\,\\
&&\llap{$\ts\frac{1}{r^2-q^{-2}}\,L\bigl(\frac{r}{q^{-2}-r^2},\frac{1}{qr}\bigr)\,L\bigl(\frac{r}{q^{-2}-r^2},\frac{1}{qr}\bigr)^t\Bigr)$ \ [using
\eqref{LLTGF}]}\\
&\ts=\frac{r^{n(n+1)}\,a^{n^2}b^n}{c^n\,(q^{-1}-qr^2)^{n(n-1)}}\,\det\Bigl(\frac{1}{r^2-q^2}\,
L\bigl(\frac{r}{q^{-2}-r^2},\frac{1}{qr}\bigr)^{-1}
L\bigl(\frac{r}{q^2-r^2},\frac{q}{r}\bigr)\Bigl(L\bigl(\frac{r}{q^{-2}-r^2},\frac{1}{qr}\bigr)^{-1}
L\bigl(\frac{r}{q^2-r^2},\frac{q}{r}\bigr)\Bigr)^t\,-\\
&&\llap{$\ts\,\frac{1}{r^2-q^{-2}}\,I\Bigr)$ \ [using \eqref{Ldet}]}\\
&\ts=\Bigl(\frac{(r^2-q^2)(r^2-q^{-2})}{q^2-q^{-2}}\Bigr)^n\det\Bigl(\frac{1}{q^{-2}-r^2}\,I\,+\,\frac{1}{r^2-q^2}\,
L\bigl(\frac{q^2-q^{-2}}{q^{-1}r-qr^{-1}},\frac{qr-(qr)^{-1}}{q^2-q^{-2}}\bigr)\,
L\bigl(\frac{q^2-q^{-2}}{q^{-1}r-qr^{-1}},\frac{qr-(qr)^{-1}}{q^2-q^{-2}}\bigr)^t\Bigr)\;\\
&&\llap{[using \eqref{wr} and \eqref{LL}]}\\
&\ts=\det\Bigl((1-\omega_+)\,I\,+\,\omega_+\,L\bigl(\frac{c}{b},\frac{a}{c}\bigr)\,
L\bigl(\frac{c}{b},\frac{a}{c}\bigr)^t\Bigr)&\llap{[using \eqref{wr} and \eqref{nuqr}]}\\
&=\det_{0\le i,j\le n-1}\ts\Bigl((1-\omega_+)\,\delta_{ij}\,+\,\omega_+\sum_{k=0}^{\min(i,j)}{i\choose k}{j\choose k}
\bigl(\frac{c}{b}\bigr)^{i+j}\,\bigl(\frac{a}{c}\bigr)^{2k}\Bigr)&\llap{[using \eqref{LLT}]}\\
&=\det_{0\le i,j\le n-1}\ts\left((1-\omega_+)\,\delta_{ij}\,+\,\omega_+\sum_{k=0}^{\min(i,j)}{i\choose k}{j\choose k}
\bigl(\frac{a}{b}\bigr)^{2k}\,\bigl(\frac{c}{b}\bigr)^{2(i-k)}\right)&\llap{[multiplying the matrix}\\[-2.4mm]
&&\llap{on the left by $V$ and right by $V^{-1}$, where $V_{ij}=\bigl(\frac{c}{b}\bigr)^i\delta_{ij}$]}.
\end{alignat*}
It can be seen from~\eqref{wr} and~\eqref{nuqr} that under the transformation $r\rightarrow r^{-1}$ and $q\rightarrow q^{-1}$,~$\frac{a}{b}$
and~$\frac{c}{b}$ are invariant, and $\omega_+$ and $\omega_-$ are interchanged.
Therefore, applying this transformation to the equality between the first and last expressions in the previous sequence
it follows that this equality remains valid if~$\omega_+$ is replaced by~$\omega_-$.
Finally, identifying $x=\bigl(\frac{a}{b}\bigr)^2$ and $y=\bigl(\frac{c}{b}\bigr)^2$ gives the result of Proposition~\ref{ASMdet},
using~\eqref{nuxy}, \eqref{MASM} and~\eqref{nuabc}.

\subsection{\protect Proof of Proposition \ref{DPPdet}}
\label{DPPdetderiv}
In this subsection, the DPP determinant formula of Proposition \ref{DPPdet} will be obtained
using the Lindstr\"{o}m--Gessel--Viennot theorem, together with a bijection between
DPPs and certain sets of nonintersecting lattice paths.  Closely related determinant formulae have
also been obtained by Mills, Robbins and Rumsey~\cite[p.~346]{MRR-ASM} and Lalonde~\cite[Thm.~3.1]{Lal-ASM}.

The bijection between DPPs and sets of nonintersecting lattice paths will be given first.
This correspondence was outlined by Lalonde~\cite[Sec.~2]{Lal-ASM} and is based on
similar such correspondences, as obtained by Gessel and Viennot~\cite{GV,GV-pre} for certain plane partitions
and Young tableaux.
Correspondences between DPPs and slightly different sets of nonintersecting lattice paths have been
given by Bressoud~\cite{Bressoud}, Krattenthaler~\cite{Kratt-DPP} and
Lalonde~\cite{Lal-ASM,Lal-DPP}, and some information about these will be given at
the end of this subsection.

Consider the directed $n\times n$ grid $G_{n,n}$ with vertices $\{(i,j)\mid i,j=0,\ldots,n-1\}$,
where $(i,j)$ is in the $i$th column from the left and $j$th row from the bottom,
and for which there are edges from $(i,j)$ to $(i+1,j)$ for each $i=0,\ldots,n-2$, $j=0,\ldots,n-1$, and from $(i,j)$ to $(i,j-1)$ for each
$i=0,\ldots,n-1$, $j=1,\ldots,n-1$.  Thus, all horizontal edges are directed rightward and all
vertical edges are directed downward.  Note that this grid differs slightly from the $n\times n$ grid used in Section~\ref{ASMdetderiv}.
In particular, the grid there was undirected, had certain additional boundary vertices and edges, and matrix-type rather than Cartesian
labeling was used.

Now let $\NILP(n)$ be the set of all sets~$P$ of nonintersecting paths on $G_{n,n}$ for which there exist
integers $0\le t\le n-1$ and $n=\lambda_0>\lambda_1>\ldots>\lambda_t>\lambda_{t+1}=0$ such that $P$~consists of
paths from $(0,\lambda_{i-1}-1)$ to $(\lambda_i,0)$ for $i=1,\ldots,t+1$.

\begin{figure}[h]
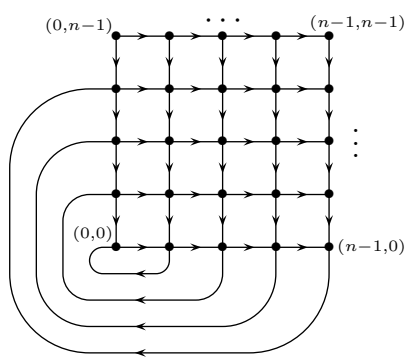
\centering
\psset{unit=3.5mm}
\pspicture(-4,-4)(9,9)
\multips(0,0)(0,2){5}{\psline[arrows=->,arrowsize=3pt,linewidth=0.5pt](0,0)(1.2,0)
\multips(1,0)(2,0){3}{\psline[arrows=->,arrowsize=3pt,linewidth=0.5pt](0,0)(2.2,0)}\psline[linewidth=0.5pt](7,0)(8,0)}
\multips(0,8)(2,0){5}{\psline[arrows=->,arrowsize=3pt,linewidth=0.5pt](0,0)(0,-1.2)
\multips(0,-1)(0,-2){3}{\psline[arrows=->,arrowsize=3pt,linewidth=0.5pt](0,0)(0,-2.2)}\psline[linewidth=0.5pt](0,-7)(0,-8)}
\psline[linewidth=0.5pt,linearc=3](8,0)(8,-4)(-4,-4)(-4,6)(0,6)
\psline[linewidth=0.5pt,linearc=2](6,0)(6,-3)(-3,-3)(-3,4)(0,4)
\psline[linewidth=0.5pt,linearc=1](4,0)(4,-2)(-2,-2)(-2,2)(0,2)
\psline[linewidth=0.5pt,linearc=0.5](2,0)(2,-1)(-1,-1)(-1,0)(0,0)
\multips(0.8,-1)(0,-1){4}{\psline[arrows=->,arrowsize=3pt,linewidth=0.5pt](0,0)(-0.1,0)}
\multirput(0,8)(2,0){5}{$\ss\bullet$}\multirput(0,6)(2,0){5}{$\ss\bullet$}\multirput(0,4)(2,0){5}{$\ss\bullet$}
\multirput(0,2)(2,0){5}{$\ss\bullet$}\multirput(0,0)(2,0){5}{$\ss\bullet$}
\rput[br](-0.1,8.1){$\sss(0,n-1)$}\rput[r](10.9,8.5){$\sss(n-1,n-1)$}
\rput[r](-0.1,0.5){$\sss(0,0)$}\rput[r](10.9,0){$\sss(n-1,0)$}
\rput[b](4.1,8.3){$\cdots$}\rput[;](9,4.3){$\vdots$}
\endpspicture
\caption{The graph obtained from $G_{n,n}$ by adding certain boundary edges.}
\label{Gnn}\end{figure}

It can be seen that each set of paths in $\NILP(n)$ corresponds to a single vertex-distinct
path from $(0,n-1)$ to $(0,0)$ on the graph shown in Figure~\ref{Gnn}, which is obtained from $G_{n,n}$ by adding certain boundary edges.

The mappings between~$\DPP(n)$ and~$\NILP(n)$,
which can be shown straightforwardly to be well-defined and bijective, are as follows.
\begin{itemize}
\item To map $P\in\NILP(n)$ to $D\in\DPP(n)$, first
number the paths of $P$ according to the height at which they start, taking the path which starts highest as path~1.
Then take part $D_{ij}$ to be $1$ plus the height of the $(j-i+1)$th rightward step of path~$i$.
\item To map $D\in\DPP(n)$ to $P\in\NILP(n)$, first let $t$ be the number of rows in~$D$ and let~$\lambda_i$ be the length
of row~$i$, as in~\eqref{DPP}.  Also define $\lambda_0=n$ and $\lambda_{t+1}=0$.  Then
obtain~$P$ by forming a path, for each $i=1,\ldots,t+1$, from $(0,\lambda_{i-1}-1)$ to
$(\lambda_i,0)$ whose~$k$th rightward step has height $D_{i,i+k-1}-1$.
\end{itemize}

It follows from these mappings that
the statistics $\nu$ and $\mu$ for DPPs, given by~\eqref{nuD} and~\eqref{muD}, can be identified with numbers of certain horizontal steps
in sets of paths as follows.
\begin{lemma}\label{numuDPP}
Let $D\in\DPP(n)$ correspond to $P\in\NILP(n)$.  Then
\begin{align*}\nu(D)&=\text{number of rightward steps in $P$ above the line }\{(i,i-1)\},\\[1mm]
\mu(D)&=\text{number of rightward steps in $P$ below the line }\{(i,i-1)\}.
\end{align*}\end{lemma}

\begin{figure}[h]
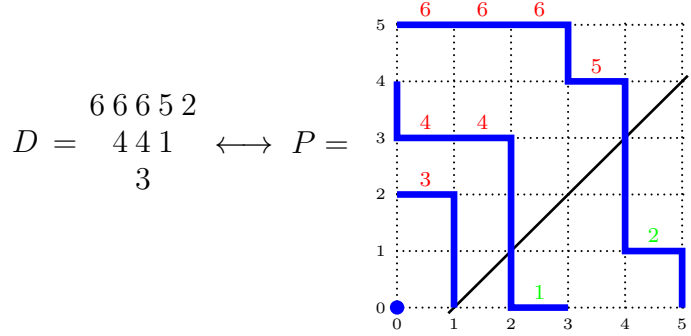
\centering
\raisebox{22mm}{$\ba{c@{\:}c@{\:}c@{\:}c@{\:}c@{\:}c}
&6&6&6&5&2\\
D\,=\,&&4&4&1\\
&&&3\ea\;\longleftrightarrow\;P=\;$}
\psset{unit=7.5mm}
\pspicture(-0.5,-0.2)(5,5.3)
\psgrid[subgriddiv=1,griddots=9,gridlabels=0](0,0)(5,5)
\psline[linewidth=1pt](0.9,-0.1)(5.1,4.1)
\psline[linewidth=2.5pt,linecolor=blue](0,5)(3,5)(3,4)(4,4)(4,1)(5,1)(5,0)
\psline[linewidth=2.5pt,linecolor=blue](0,4)(0,3)(2,3)(2,0)(3,0)
\psline[linewidth=2.5pt,linecolor=blue](0,2)(1,2)(1,0)
\psdots[dotscale=1.5,linecolor=blue](0,0)
\rput[b](0.5,5.15){$\color{red}\ss6$}\rput[b](1.5,5.15){$\color{red}\ss6$}\rput[b](2.5,5.15){$\color{red}\ss6$}\rput[b](3.5,4.15){$\color{red}\ss5$}
\rput[b](4.5,1.15){$\color{green}\ss2$}
\rput[b](0.5,3.15){$\color{red}\ss4$}\rput[b](1.5,3.15){$\color{red}\ss4$}
\rput[b](0.5,2.15){$\color{red}\ss3$}\rput[b](2.5,0.15){$\color{green}\ss1$}
\rput[t](0,-0.2){$\sss0$}\rput[t](1,-0.2){$\sss1$}\rput[t](2,-0.2){$\sss2$}\rput[t](3,-0.2){$\sss3$}\rput[t](4,-0.2){$\sss4$}\rput[t](5,-0.2){$\sss5$}
\rput[r](-0.2,0){$\sss0$}\rput[r](-0.2,1){$\sss1$}\rput[r](-0.2,2){$\sss2$}\rput[r](-0.2,3){$\sss3$}\rput[r](-0.2,4){$\sss4$}\rput[r](-0.2,5){$\sss5$}
\endpspicture
\caption{A DPP and corresponding set of nonintersecting lattice paths.}
\label{DPPPEx}\end{figure}

The set of nonintersecting lattice paths which corresponds to the example of~\eqref{DPPEx} (with $n=6$)
is shown in Figure~\ref{DPPPEx}, where each part of $D$ is shown above its corresponding rightward step, with the
nonspecial and special parts in red and green respectively, and the line $\{(i,i-1)\}$ also shown.

As a further example, mapping $\DPP(3)$, as given in \eqref{DPP3}, to $\NILP(3)$ gives
\psset{unit=5mm}
\begin{equation}\label{NILP3}\NILP(3)\,=\,\left\{\raisebox{-4mm}{
\pspicture(0,0)(2.8,2)\psgrid[subgriddiv=1,griddots=6,gridlabels=0](0,0)(2,2)
\psline[linewidth=1.7pt,linecolor=blue](0,2)(0,0)\rput(2.3,0.5){,}\endpspicture
\pspicture(0,0)(2.8,2)\psgrid[subgriddiv=1,griddots=6,gridlabels=0](0,0)(2,2)
\psline[linewidth=1.7pt,linecolor=blue](0,2)(2,2)(2,0)\psline[linewidth=1.7pt,linecolor=blue](0,1)(1,1)(1,0)
\psdots[dotscale=0.9,linecolor=blue](0,0)\rput(2.3,0.5){,}\endpspicture
\pspicture(0,0)(2.8,2)\psgrid[subgriddiv=1,griddots=6,gridlabels=0](0,0)(2,2)
\psline[linewidth=1.7pt,linecolor=blue](0,2)(0,1)(1,1)(1,0)
\psdots[dotscale=0.9,linecolor=blue](0,0)\rput(2.3,0.5){,}\endpspicture
\pspicture(0,0)(2.8,2)\psgrid[subgriddiv=1,griddots=6,gridlabels=0](0,0)(2,2)
\psline[linewidth=1.7pt,linecolor=blue](0,2)(2,2)(2,0)\psline[linewidth=1.7pt,linecolor=blue](0,1)(0,0)\rput(2.3,0.5){,}\endpspicture
\pspicture(0,0)(2.8,2)\psgrid[subgriddiv=1,griddots=6,gridlabels=0](0,0)(2,2)
\psline[linewidth=1.7pt,linecolor=blue](0,2)(1,2)(1,0)\psdots[dotscale=0.9,linecolor=blue](0,0)\rput(2.3,0.5){,}\endpspicture
\pspicture(0,0)(2.8,2)\psgrid[subgriddiv=1,griddots=6,gridlabels=0](0,0)(2,2)
\psline[linewidth=1.7pt,linecolor=blue](0,2)(1,2)(1,1)(2,1)(2,0)
\psline[linewidth=1.7pt,linecolor=blue](0,1)(0,0)\rput(2.3,0.5){,}\endpspicture
\pspicture(0,0)(2.2,2)\psgrid[subgriddiv=1,griddots=6,gridlabels=0](0,0)(2,2)
\psline[linewidth=1.7pt,linecolor=blue](0,2)(1,2)(1,0)(2,0)
\psline[linewidth=1.7pt,linecolor=blue](0,1)(0,0)\endpspicture}
\right\},\end{equation}
where the elements are listed in the order corresponding to that used in~\eqref{DPP3}.

Now consider an acyclic directed graph~$\mathcal{G}$.
Let a weight be assigned to each edge of~$\mathcal{G}$, and
define the weight~$W(p)$ of a path~$p$ on~$\mathcal{G}$ to be the product of the weights of its edges.
For vertices~$u$ and~$v$ of~$\mathcal{G}$, let $\paths(u,v)$ denote the set of all
paths on~$\mathcal{G}$ from~$u$ to~$v$.   For vertices $u_1,\ldots,u_n,v_1,\ldots,v_n$ of~$\mathcal{G}$,
let $\mathcal{N}_{\mathcal{G}}(u_1,\ldots,u_n;v_1,\ldots,v_n)$ denote the set of all sets~$P$ of
paths on~$\mathcal{G}$ such that $P$ consists of a path of $\paths(u_i,v_i)$ for each $i=1,\ldots,n$, and
different paths of~$P$ are nonintersecting.
The Lindstr\"{o}m--Gessel--Viennot theorem~\cite{GV,GV-pre,Lind} can now be stated as follows.
\begin{theorem*}[Lindstr\"{o}m, Gessel, Viennot]
If $\mathcal{N}_{\mathcal{G}}(u_{\sigma_1},\ldots,u_{\sigma_n};v_1,\ldots,v_n)$ is empty for each permutation $\sigma\in\mathcal{S}_n$
other than the identity, then
\begin{equation}\label{LGV}\sum_{P\in\mathcal{N}_{\mathcal{G}}(u_1,\ldots,u_n;v_1,\ldots,v_n)}\,\prod_{p\in P}W(p)\,=\det_{1\le i,j\le n}\Biggl(
\sum_{p\in\paths(u_j,v_i)}W(p)\Biggr).\end{equation}
\end{theorem*}
Briefly, this can be proved by writing the RHS as
$\sum_{\sigma\in\mathcal{S}_n}\sum_{p_1\in\paths(u_{\sigma_1},v_1)}\ldots\sum_{p_n\in\paths(u_{\sigma_n},v_n)}$
$\mathrm{sgn}(\sigma)W(p_1)\ldots W(p_n)$.  A certain
involution~$\phi$ is then constructed on the set
$\{(p_1,\ldots,p_n)\in\paths(u_{\sigma_1},v_1)\times\ldots\times\paths(u_{\sigma_n},v_n)\mid\sigma\in\mathcal{S}_n$,
there exist $i$ and $j$ for which $p_i$ and $p_j$ intersect$\}$.
This involution has the properties that
if $\phi(p_1,\ldots,p_n)=(p'_1,\ldots,p'_n)$, with
$p_i\in\paths(u_{\sigma_i},v_i)$ and $p'_i\in\paths(u_{\sigma'_i},v_i)$ for each $i=1,\ldots,n$, then
$\mathrm{sgn}(\sigma)=-\mathrm{sgn}(\sigma')$ and $W(p_1)\ldots W(p_n)=W(p'_1)\ldots W(p'_n)$.  It then
follows that each corresponding pair of terms on the RHS of~\eqref{LGV} cancels, leaving only the terms which give the LHS.
For further details of this proof see Gessel and Viennot~\cite[Sec.~2]{GV-pre} or Lindstr\"{o}m~\cite[Lemma~1]{Lind}.
For an account of the origins of the theorem, or closely related results, see for example Krattenthaler~\cite[Footnote~10]{Kratt-det2}.
For an account of the relationship with free fermionic methods see for example Zinn-Justin~\cite[Sec.~1]{hdr}.

The set $\NILP(n)$ can be written as
\begin{multline}\label{NILP}
\NILP(n)=\\
\bigcup_{n-1\ge\lambda_1>\ldots>\lambda_t\ge1}\!
\mathcal{N}_{G_{n,n}}\Bigl((0,n-1),(0,\lambda_1-1),\ldots,(0,\lambda_t-1);(\lambda_1,0)\ldots,(\lambda_t,0),(0,0)\Bigr).
\end{multline}
Now assign the horizontal edge in $G_{n,n}$ from $(i,j)$ to $(i+1,j)$ a weight of~$x$ for $i\le j$ and
a weight of~$y$ for $i>j$, and assign each vertical edge a weight of~1.
Then the bijection between $\DPP(n)$ and $\NILP(n)$, together with
\eqref{ZDPP}, \eqref{LGV} and Lemma~\ref{numuDPP} (and the fact that the condition for the
validity of~\eqref{LGV} is satisfied), gives
\begin{equation}\label{DPPLGV1}
Z_\DPP(n,x,y,1)=
\sum_{n-1\ge\lambda_1>\ldots>\lambda_t\ge1}\;
\det_{\rule{0ex}{2.3ex}\substack{i=0,\lambda_t,\ldots,\lambda_1\\j=\lambda_t-1,\ldots,\lambda_1-1,n-1}}
\Biggl(\sum_{p\in\paths((0,j),(i,0))}W(p)\Biggr).
\end{equation}
Let $S$ be the $n\times n$  subdiagonal matrix,
\begin{equation}\label{S}S_{ij}=\delta_{i,j+1}.\end{equation}
Then, for any matrix $A$ with rows and columns indexed by $\{0,\ldots,n-1\}$,
\begin{equation}\label{detidentity}\det(A-S)=\sum_{T\subset\{1,\ldots,n-1\}}\det A_{\{0\}\cup T,(T-1)\cup\{n-1\}},\end{equation}
where  $A_{\{0\}\cup T,(T-1)\cup\{n-1\}}$ is the
submatrix of $A$ formed by restricting the rows and columns to those indexed by
$\{0\}\cup T$ and $\{t-1\mid t\in T\}\cup\{n-1\}$ respectively.
Briefly, this identity can be obtained by writing $\det(A-S)=
\sum_{\sigma\in\mathcal{S}_n}\mathrm{sgn}(\sigma)A_{0,\sigma_0}\prod_{i=1}^{n-1}(A_{i,\sigma_i}-\delta_{i,\sigma_i+1})=
\sum_{T\subset\{1,\ldots,n-1\}}(-1)^{n+|T|+1}\sum_{\sigma\in\mathcal{S}_n}\mathrm{sgn}(\sigma)A_{0,\sigma_0}\prod_{i\in T}A_{i,\sigma_i}
\prod_{i\in\{1,\ldots,n-1\}\setminus T}\delta_{i,\sigma_i+1}$.  The result then follows using
a certain decomposition property of $\mathrm{sgn}(\sigma)$.

Applying~\eqref{detidentity} to~\eqref{DPPLGV1} now gives
\begin{equation}\label{DPPLGV2}
Z_\DPP(n,x,y,1)=
\det_{0\le i,j\le n-1}\Biggl(-\delta_{i,j+1}+\sum_{p\in\paths((0,j),(i,0))}W(p)\Biggr).
\end{equation}
The paths of $\paths((0,j),(i,0))$ can be obtained by joining any path of $\paths((0,j),(k,k-1))$ to any path
of $\paths((k,k-1),(i,0))$,
for any $0\le k\le\min(i,j+1)$.  Each path of~$\paths((0,j),(k,k-1))$ has $k$ horizontal edges, each
with weight~$x$, and $j-k+1$ vertical edges, each with weight~1, so there are $j+1\choose k$ such paths, each with weight~$x^k$.
Similarly, $\paths((k,k-1),(i,0))$ contains $i-1\choose i-k$ paths, each with weight~$y^{i-k}$.
This is shown diagrammatically in Figure~\ref{DPPdetfig}.
(Note that the case $k=0$ is slightly exceptional: an additional edge from $(0,0)$ to $(0,-1)$ is then considered, with
$\paths((0,j),(0,-1))$ containing a single path, and
$\paths((0,-1),(i,0))$ containing a single path for $i=0$ and empty otherwise.)
It now follows that
\begin{equation}\label{DPPwp}
\sum_{p\in\paths((0,j),(i,0))}W(p)=\sum_{k=0}^{\min(i,j+1)}{i-1\choose i-k}{j+1\choose k}x^k\,y^{i-k}.
\end{equation}

\begin{figure}[h]\centering
\psset{unit=4mm}
\pspicture(-3.5,-1.1)(10,10)
\psgrid[subgriddiv=1,griddots=6,gridlabels=0](0,0)(10,10)
\psline[linewidth=1pt](0.9,-0.1)(10.1,9.1)
\rput(7,-0.9){$\sss i-k$}\rput(2.5,-0.9){$\sss k$}\rput(4.5,-1.8){$\sss i$}
\rput(-1.5,2){$\sss k-1$}\rput(-1.5,6){$\sss j-k+1$}\rput(-3.2,4){$\sss j$}
\psline[linewidth=1pt]{<-}(0,-0.9)(2.1,-0.9)\psline[linewidth=1pt]{->}(2.9,-0.9)(4.9,-0.9)
\psline[linewidth=1pt]{<-}(5.1,-0.9)(6.1,-0.9)\psline[linewidth=1pt]{->}(7.9,-0.9)(9,-0.9)
\psline[linewidth=1pt]{<-}(0,-1.8)(4.1,-1.8)\psline[linewidth=1pt]{->}(4.9,-1.8)(9,-1.8)
\psline[linewidth=1pt]{<-}(-1.5,0)(-1.5,1.5)\psline[linewidth=1pt]{->}(-1.5,2.5)(-1.5,3.9)
\psline[linewidth=1pt]{<-}(-1.5,4.1)(-1.5,5.5)\psline[linewidth=1pt]{->}(-1.5,6.5)(-1.5,8)
\psline[linewidth=1pt]{<-}(-3.2,0)(-3.2,3.5)\psline[linewidth=1pt]{->}(-3.2,4.5)(-3.2,8)
\psline[linewidth=1.7pt,linecolor=blue](0,8)(2,8)(2,7)(4,7)(4,4)(5,4)(5,3)(7,3)(7,2)(8,2)(8,0)(9,0)
\rput[b](0.5,8.2){$\sss x$}\rput[b](1.5,8.2){$\sss x$}\rput[b](2.5,7.2){$\sss x$}\rput[b](3.5,7.2){$\sss x$}\rput[b](4.5,4.2){$\sss x$}
\rput[b](5.5,3.2){$\sss y$}\rput[b](6.5,3.2){$\sss y$}\rput[b](7.5,2.2){$\sss y$}\rput[b](8.5,0.2){$\sss y$}
\endpspicture
\caption{Derivation of~\eqref{DPPwp}.}\label{DPPdetfig}\end{figure}
Finally, combining~\eqref{DPPLGV2} and~\eqref{DPPwp} gives
the result of Proposition~\ref{DPPdet}.

An alternative set $\NILP'(n)$ of sets of nonintersecting lattice paths which is also in bijection with $\DPP(n)$ will now
be described.  This is not needed for the main results of the paper, but is included for completeness, and since
it matches certain nonintersecting lattice paths in the literature, in particular those
of Bressoud~\cite[p.~108 \& Exercise~3.4.9]{Bressoud} and Krattenthaler~\cite{Kratt-DPP}.
Some differences between $\NILP(n)$ and $\NILP'(n)$
are that they use an $n\times n$ and an $(n-1)\times(n+1)$ grid respectively,
and that the endpoints of the paths are closely related to the lengths of rows of DPPs for $\NILP(n)$
and to the first entries of rows of DPPs for $\NILP'(n)$.

Let $G_{n-1,n+1}$ be the directed $(n-1)\times(n+1)$ grid with vertices $\{(i,j)\mid i=1,\ldots,n-1,\ j=-1,\ldots,n-1\}$,
and in which each horizontal edge is directed rightward and each vertical edge is directed downward.
Let $\NILP'(n)$ be the set of all sets~$P'$ of nonintersecting paths on~$G_{n-1,n+1}$ for which
there exist integers $0\le t\le n-1$ and $n-1\ge\delta_1>\ldots>\delta_t\ge1$ such that $P'$~consists of
paths from $(1,\delta_i)$ to $(\delta_i,-1)$ for $i=1,\ldots,t$.

To map $P'\in\NILP'(n)$ to $D\in\DPP(n)$, first number the paths of $P'$ according to the height
at which they start, taking the path which starts highest as path~1. Then take $D_{ii}$ to be $1$ plus the initial height of
path $i$, and $D_{ij}$ for $j>i$ to be $1$ plus the height of the $(j-i)$th rightward step of path~$i$,
excluding steps at height $-1$.
To map $D\in\DPP(n)$ to $P'\in\NILP'(n)$, let $t$ be the number of rows in~$D$ and
obtain~$P'$ by forming a path, for each $i=1,\ldots,t$, from $(1,D_{ii}-1)$ to $(D_{ii}-1,-1)$
whose~$k$th rightward step has height $D_{i,i+k}-1$, adding further steps at height $-1$ if necessary.
It can be seen that
$\nu(D)$ is the number of paths in $P'$ plus the number of rightward steps in $P'$ above
the line $\{(i,i-1)\}$,
and that~$\mu(D)$ is the number of rightward steps in $P'$ below the line $\{(i,i-1)\}$
and with nonnegative height.
The alternative set of nonintersecting lattice paths which corresponds to the example
of~\eqref{DPPEx} (with $n=6$) is shown in Figure~\ref{DPPAltEx}.
\begin{figure}[h]\centering
\raisebox{22mm}{$\ba{c@{\:}c@{\:}c@{\:}c@{\:}c@{\:}c}
&6&6&6&5&2\\
D\,=\,&&4&4&1\\
&&&3\ea\;\longleftrightarrow\;P'=\;$}
\psset{unit=7.5mm}
\pspicture(0.5,-1.2)(5,5.3)
\psgrid[subgriddiv=1,griddots=9,gridlabels=0](1,-1)(5,5)
\psline[linewidth=1pt](0.9,-0.1)(5.1,4.1)
\psline[linewidth=2.5pt,linecolor=blue](1,5)(3,5)(3,4)(4,4)(4,1)(5,1)(5,-1)
\psline[linewidth=2.5pt,linecolor=blue](1,3)(2,3)(2,0)(3,0)(3,-1)
\psline[linewidth=2.5pt,linecolor=blue](1,2)(1,2)(1,-1)(2,-1)
\rput[b](0.9,5.15){$\color{red}\ss6$}\rput[b](1.5,5.15){$\color{red}\ss6$}\rput[b](2.5,5.15){$\color{red}\ss6$}\rput[b](3.5,4.15){$\color{red}\ss5$}
\rput[b](4.5,1.15){$\color{green}\ss2$}
\rput[b](0.83,3.15){$\color{red}\ss4$}\rput[b](1.5,3.15){$\color{red}\ss4$}
\rput[b](0.83,2.15){$\color{red}\ss3$}\rput[b](2.5,0.15){$\color{green}\ss1$}
\rput[t](1,-1.2){$\sss1$}\rput[t](2,-1.2){$\sss2$}\rput[t](3,-1.2){$\sss3$}\rput[t](4,-1.2){$\sss4$}\rput[t](5,-1.2){$\sss5$}
\rput[r](0.7,-1){$\sss-1$}\rput[r](0.7,0){$\sss0$}\rput[r](0.7,1){$\sss1$}
\rput[r](0.7,2){$\sss2$}\rput[r](0.7,3){$\sss3$}\rput[r](0.7,4){$\sss4$}\rput[r](0.7,5){$\sss5$}
\endpspicture
\caption{A DPP and corresponding alternative set of nonintersecting lattice paths.}
\label{DPPAltEx}\end{figure}
Now assign the horizontal edge from $(i,j)$ to $(i+1,j)$ a weight of~$x$ for $i\le j$,
a weight of~$y$ for $i>j\ge0$ and a weight of 1 for $j=-1$, and assign each vertical edge a weight of~1.
Then the Lindstr\"{o}m--Gessel--Viennot theorem~\eqref{LGV}, together with the identity
$\det(I+A)=\sum_{T\subset K}\det A_{T,T}$
(where $A$ is a matrix with rows and columns indexed by $K$ and
$A_{T,T}$ is the submatrix formed by restricting the rows and columns to those indexed by $T$), gives
$Z_\DPP(n,x,y,1)=\det M'_\DPP(n,x,y,1)$, with
\begin{equation}\label{Md}
M'_\DPP(n,x,y,1)_{ij}=\delta_{ij}+\sum_{k=0}^{i-1}\,\sum_{l=0}^{\min(j,k)}{j\choose l}{k\choose l}x^{l+1}\,y^{k-l},\end{equation}
where either $1\le i,j\le n-1$ or $0\le i,j\le n-1$ can be used.
Note that $M'_\DPP(n,1,1,1)_{ij}=\delta_{ij}+{i+j\choose i-1}$, and that
$(I-S)M'_\DPP(n,x,y,1)=M_\DPP(n,x,y,1)(I-S^t)$,
with~$S$ given by~\eqref{S} and $M'_\DPP(n,x,y,1)_{ij}$ taken with $0\le i,j\le n-1$.
The formula $Z_\DPP(n,1,1,1)=\det M'_\DPP(n,1,1,1)$
also follows from a result of Andrews~\cite[Thm.~3]{And-weakMac},
and the formula $Z_\DPP(n,1,y,1)=\det M'_\DPP(n,1,y,1)$
follows from results of Mills, Robbins and Rumsey~\cite[pp.~50 \&~54]{MRR-symPP}.

Lalonde has shown that $\DPP(n)$ is also in bijection with certain sets
of~$n$ lattice paths, termed top-bottom or TB configurations,
on an equilateral triangular lattice of side length~$n$. See~\cite[Sec.~3]{Lal-ASM} and~\cite[Sec.~2]{Lal-DPP} for details.
The TB configuration which corresponds to $P\in\NILP(n)$
is obtained by applying a certain transformation to the segments of paths of $P$ below the line
$\{(i,i-1)\}$.  By using TB configurations and a Lindstr\"{o}m--Gessel--Viennot-type theorem, Lalonde~\cite[Thm.~3.1]{Lal-ASM}
showed that $Z_\DPP(n,x,y,1)=\det M''_\DPP(n,x,y,1)$, where
\begin{equation}\label{Mdd}M''_\DPP(n,x,y,1)_{ij}={j+1\choose i}x^i-{i-1\choose i-j-1}(-y)^{i-j-1},\end{equation}
with $0\le i,j\le n-1$.  Note that
$B(y)M''_\DPP(n,x,y,1)=M_\DPP(n,x,y,1)$,
where $B(y)_{ij}={i-1\choose i-j}y^{i-j}$.

\begin{figure}[h]
$\vcenter{
\hbox{$\ba{c@{\:}c@{\:}c@{\:}c@{\:}c@{\:}c}
&6&6&6&5&2\\
&&4&4&1\\
&&&3\ea\leftrightarrow\;\;$}}
\vcenter{
\psset{unit=5.2mm}
\hbox{\pspicture(0.5,-1.2)(5,5.3)
\psgrid[subgriddiv=1,griddots=9,gridlabels=0](1,-1)(5,5)
\psline[linewidth=2.5pt,linecolor=blue](1,5)(3,5)(3,4)(4,4)(4,1)(5,1)(5,-1)
\psline[linewidth=2.5pt,linecolor=blue](1,3)(2,3)(2,0)(3,0)(3,-1)
\psline[linewidth=2.5pt,linecolor=blue](1,2)(1,2)(1,-1)(2,-1)
\rput[b](0.9,5.2){$\ss6$}\rput[b](1.5,5.2){$\ss6$}\rput[b](2.5,5.2){$\ss6$}\rput[b](3.5,4.2){$\ss5$}
\rput[b](4.5,1.2){$\ss2$}
\rput[b](0.8,3.2){$\ss4$}\rput[b](1.5,3.2){$\ss4$}
\rput[b](0.8,2.2){$\ss3$}\rput[b](2.5,0.2){$\ss1$}
\endpspicture}}
\;\;\leftrightarrow\;\;
\vcenter{\hbox{\includegraphics[scale=0.4]{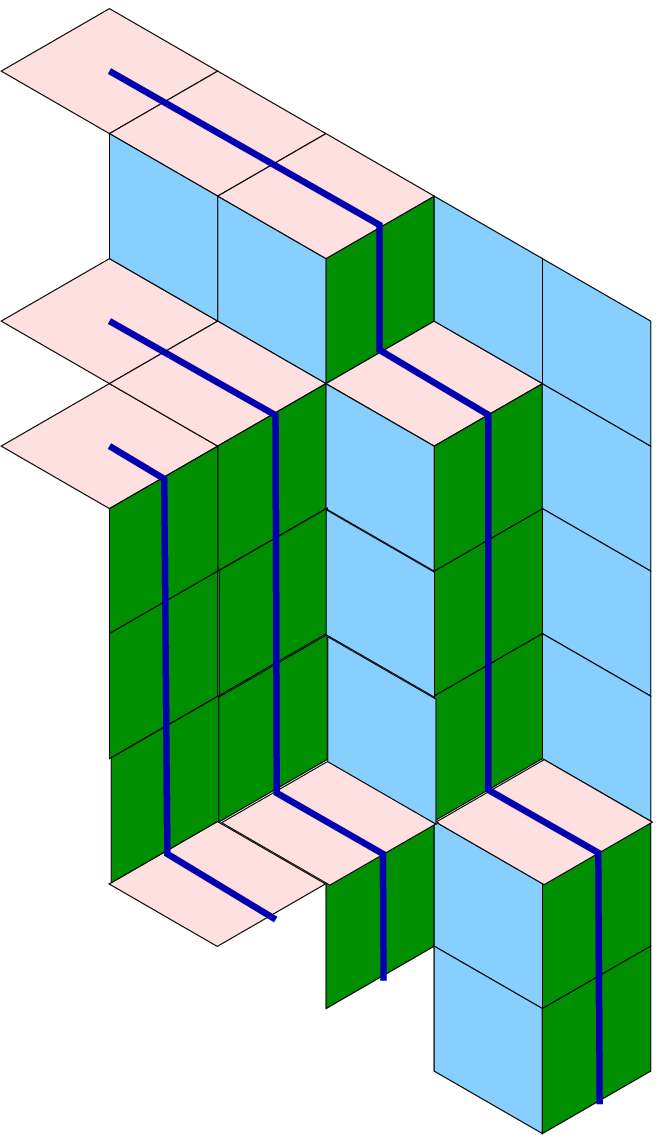}}}
\;\;\leftrightarrow\;\;
\vcenter{\hbox{\includegraphics[scale=0.4]{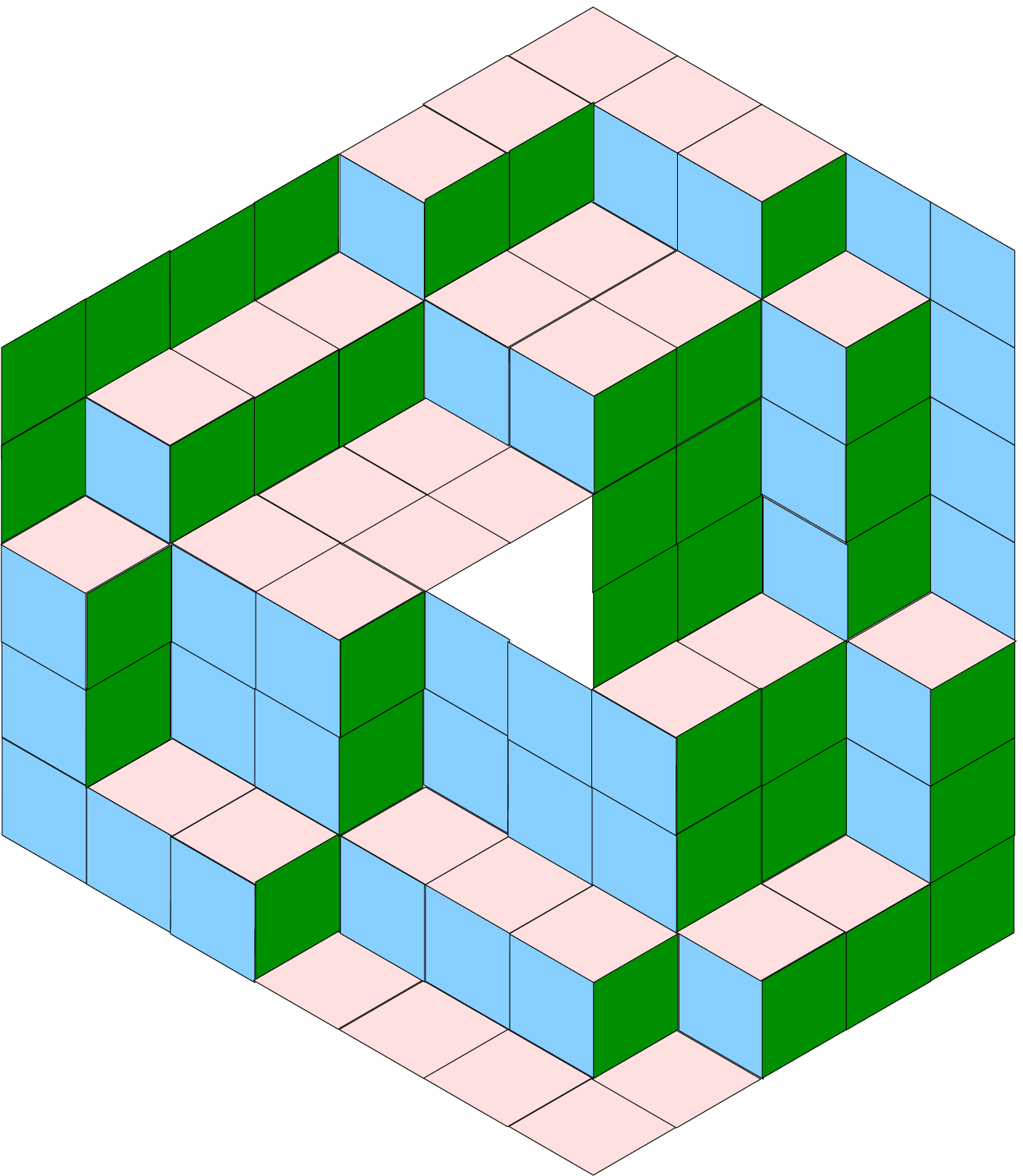}}}
$
\caption{A DPP and corresponding cyclically symmetric rhombus tiling of a punctured hexagon.}
\label{hex}
\end{figure}

Finally, Krattenthaler~\cite{Kratt-DPP} has shown that $\DPP(n)$ is also in bijection with the set of cyclically symmetric rhombus tilings of a
hexagon with alternating sides of length $n-1$ and~$n+1$, from which a central equilateral triangular hole of side length~2 has been removed.
This bijection is best understood in terms of the alternative
set $\NILP'(n)$ of sets of nonintersecting lattice paths. The
correspondence between the example of~\eqref{DPPEx} (with $n=6$) and a rhombus tiling of a punctured hexagon is shown in Figure~\ref{hex}.

\subsection{\protect Proof of Proposition \ref{ASMdetDPPdet}}\label{ASMdetDPPdetproof}
In this subsection, the equality of the ASM and DPP determinants in Proposition~\ref{ASMdetDPPdet} will be proved by
obtaining a certain relation between the generating functions associated with the entries of the
respective matrices.  This will then imply a relation between the matrices themselves,
and the equality of their determinants.

Define the generating functions
\begin{align}g_\ASM(x,y;u,v)&=\frac{1-\omega}{1-uv}+\frac{\omega}{1-yu-v-(x-y)uv},\\[2mm]
g_\DPP(x,y;u,v)&=\frac{1-u}{(1-v)(1-uv)}+\frac{xu}{(1-v)(1-yu-v-(x-y)uv)}.\end{align}
Then
\begin{align}\label{gASM}M_\ASM(n,x,y,1)_{ij}&=[u^iv^j]g_\ASM(x,y;u,v),\\[2mm]
\label{gDPP}M_\DPP(n,x,y,1)_{ij}&=[u^iv^j]g_\DPP(x,y;u,v).\end{align}
Using the equation~\eqref{nuxy} satisfied by $\omega$, it is seen that
\begin{equation}\label{fASMMfDPP}(1+(x-\omega y-1)u)g_\ASM(x,y;u,v)=(1+(\omega-1)v)g_\DPP(x,y;u,v).\end{equation}
(In particular, the difference between the two sides of~\eqref{fASMMfDPP} contains
the LHS of~\eqref{nuxy} as a factor, and therefore vanishes.)
It follows that
\begin{equation}\label{MASMMDPP}
(I+(x-\omega y-1)S)\,M_\ASM(n,x,y,1)=M_\DPP(n,x,y,1)\,(I+(\omega-1)S^t),\end{equation}
where $S$ is given by~\eqref{S}, and this now implies the result of Proposition~\ref{ASMdetDPPdet}.

\section{The refined enumeration ($z$ arbitrary)}\label{refenum}
In this section, Theorem~\ref{MRRC1} will be proved for $z$ arbitrary in~\eqref{MRRC2}.
The same approach will be taken as in Section~\ref{unref}, i.e., formulae for the ASM and DPP generating functions~\eqref{ZASM}
and~\eqref{ZDPP} as (up to a common normalization factor) determinants of $n\times n$ matrices
will be obtained, and it will then be shown that the determinants are equal.

Let $\omega$ again be a solution of~\eqref{nuxy}, and define $n\times n$ matrices
\begin{multline}\label{MASMref}
M_\ASM(n,x,y,z)_{ij}=\\[-4mm]
(1-\omega)\,\delta_{ij}+\omega
\begin{cases}
\ds\sum_{k=0}^{\min(i,j)}{i\choose k}{j\choose k}x^k\,y^{i-k},& j\le n-2\\[5mm]
\ds\sum_{k=0}^i\,\sum_{l=0}^k{i\choose k}{n-l-2\choose k-l}x^k\,y^{i-k}\,z^{l+1},& j=n-1,
\end{cases}
\end{multline}
\begin{multline}
\label{MDPPref}M_\DPP(n,x,y,z)_{ij}=\\
-\delta_{i,j+1}+
\begin{cases}
\ds\sum_{k=0}^{\min(i,j+1)}{i-1\choose i-k}{j+1\choose k}x^k\,y^{i-k},&j\le n-2\\
\ds(1+\omega(z-1))\sum_{k=0}^i\,\sum_{l=0}^k{i-1\choose i-k}{n-l-1\choose k-l}x^k\,y^{i-k}\,z^l,&j=n-1,
\end{cases}
\end{multline}
with $0\le i,j\le n-1$. Note that $\omega$ depends on~$x$ and~$y$, and that
\eqref{MASMref} and \eqref{MDPPref} reduce to~\eqref{MASM} and~\eqref{MDPP} for $z=1$.

Formulae for the ASM and DPP generating functions in terms of determinants,
and the equality of these determinants, can now be stated as follows.
\begin{proposition}\label{ASMdetref}
Let $Z_\ASM(n,x,y,z)$ and $M_\ASM(n,x,y,z)$ be given by \eqref{ZASM} and \eqref{MASMref}.  Then
\[(1+\omega(z-1))Z_\ASM(n,x,y,z)=\det M_\ASM(n,x,y,z).\]
\end{proposition}
\begin{proposition}\label{DPPdetref}
Let $Z_\DPP(n,x,y,z)$ and $M_\DPP(n,x,y,z)$ be given by \eqref{ZDPP} and \eqref{MDPPref}.  Then
\[(1+\omega(z-1))Z_\DPP(n,x,y,z)=\det M_\DPP(n,x,y,z).\]
\end{proposition}
\begin{proposition}\label{ASMdetDPPdetref}
Let $M_\ASM(n,x,y,z)$ and $M_\DPP(n,x,y,z)$ be given by~(\ref{MASMref}) and~(\ref{MDPPref}). Then
\[\det M_\ASM(n,x,y,z)=\det M_\DPP(n,x,y,z).\]
\end{proposition}
(In fact, Proposition~\ref{DPPdetref} is valid for arbitrary $\omega$, whereas the validity of
Propositions~\ref{ASMdetref} and~\ref{ASMdetDPPdetref} depends on $\omega$ being a solution of~\eqref{nuxy}.)

The validity of Theorem~\ref{MRRC1} for $z$ arbitrary in~\eqref{MRRC2} follows immediately from Propositions~\ref{ASMdetref}--\ref{ASMdetDPPdetref}.
The proofs of these propositions will be given in the ensuing three subsections.

\subsection{\protect Proof of Proposition \ref{ASMdetref}}\label{ASMdetderivref}
In this subsection, the ASM determinant formula of Proposition \ref{ASMdetref} will be obtained using
the same approach as in Section~\ref{ASMdetderiv},
i.e., the Izergin--Korepin formula~\eqref{Idet1} and the bijection between ASMs and
certain six-vertex model configurations will be used.

Consider again the set $\SVDWBC(n)$, as defined in Section~\ref{ASMdetderiv}. For any $C\in\SVDWBC(n)$,
denote the numbers of vertex configurations of~$C$ in row~1 of the grid
of types~$a_1$,~$b_1$ and $c_1$ as~$\widetilde{N}_a(C)$,~$\widetilde{N}_b(C)$ and $\widetilde{N}_c(C)$ respectively.
It can be seen that these numbers satisfy
\begin{equation}\label{paramref}\widetilde{N}_a(C)+\widetilde{N}_b(C)=n-1,\qquad \widetilde{N}_c(C)=1.\end{equation}
Also, it follows from the bijection between $\ASM(n)$ and $\SVDWBC(n)$ that the statistic~$\rho$ for ASMs, as given by~\eqref{rhoA},
satisfies the following result.
\begin{lemma}\label{rhoASM}
Let $A\in\ASM(n)$ correspond to $C\in\SVDWBC(n)$.  Then
\[\rho(A)=\widetilde{N}_a(C).\]
\end{lemma}
Now let weights $a$, $b$ and $c$ be again defined by~\eqref{wr},
for parameters $r$ and~$q$.  Also, using~\eqref{wuv} and a further parameter $s$, define weights $\tilde{a}$, $\tilde{b}$ and $\tilde{c}$ by
\begin{equation}\label{wrs}\ds \tilde{a}=\bar{a}(s,r)=sq-\frac{1}{rq},\qquad \tilde{b}=\bar{b}(s,r)=\frac{s}{q}-\frac{q}{r},\qquad
\tilde{c}=\bar{c}(s,r)=\Bigl(q^2-\frac{1}{q^2}\Bigr)
\sqrt{\frac{s}{r}}.\end{equation}

The partition function \eqref{ZSVDWBC} with these weights and $u_1=s$, $u_2=\cdots=u_n=v_1=\cdots=v_n=r$ is,
using~\eqref{Na1a2}--\eqref{Nc1c2}, given by
\begin{equation}Z_\SVDWBC(s,r,\ldots,r;r,\ldots,r)=\!\sum_{C\in\SVDWBC(n)}\!\ts
a^{2N_a(C)}\,\bigl(\frac{\tilde{a}}{a}\bigr)^{\widetilde{N}_a(C)}\,
b^{2N_b(C)}\,\bigl(\frac{\tilde{b}}{b}\bigr)^{\widetilde{N}_b(C)}\,c^{2N_c(C)+n-1}\:\tilde{c}.\end{equation}
The bijection between $\ASM(n)$ and $\SVDWBC(n)$, together with~\eqref{ZASM}, \eqref{paramsum}, \eqref{paramref} and
Lemmas~\ref{numuASM} and~\ref{rhoASM}, now gives
\begin{equation}\label{ZCZASMref}\ts Z_\SVDWBC(s,r,\ldots,r;r,\ldots,r)=b^{(n-1)^2}\,\tilde{b}^{n-1}\,c^{n-1}\,\tilde{c}\;
Z_\ASM\bigl(n,\bigl(\frac{a}{b}\bigr)^2,\bigl(\frac{c}{b}\bigr)^2,\frac{\tilde{a}\,b}{\rule{0ex}{1.65ex}a\,\tilde{b}}\bigr).\end{equation}
The following variation of \eqref{DDdet} holds for a power series $f(u,v)$ and parameters $u_0,\ldots,u_{n-1}$, $v_0,\ldots,v_{n-1}$:
\begin{multline}\label{DDdetref}
\frac{\ds\det_{0\le i,j\le n-1}\Bigl(f(u_i,v_j)\Bigr)}{\prod_{1\le i<j\le n-1}(u_j-u_i)\,\prod_{0\le i<j\le n-1}(v_j-v_i)}=\\
\det_{0\le i,j\le n-1}\left(\begin{cases}f[u_0][v_0,\ldots,v_j],&i=0\\
f[u_1,\ldots,u_i][v_0,\ldots,v_j],&i\ge1\end{cases}\right),\end{multline}
where the divided differences are again given by~\eqref{DD1}.
Now define, for parameters $\alpha$, $\beta$, $\gamma$ and $\delta$, the $n\times n$ lower triangular matrix
\begin{align}K(\alpha,\beta,\gamma,\delta)_{ij}&=\begin{cases}{i\choose j}\,\alpha^i\,\beta^j,&i\le n-2\\[2mm]
\gamma\,\delta^j,&i=n-1\end{cases}\\
\notag&=\begin{cases}[u^iv^j]\,\frac{1}{1-\alpha u(1+\beta v)},&i\leq n-2\\
[v^j]\,\frac{\gamma}{1-\delta v},&i=n-1,\end{cases}\end{align}
with $0\le i,j\le n-1$.

Then the following properties are satisfied:
\begin{align}\label{LK}\hspace*{10mm}\bigl(L(\alpha,\beta)\,K(\alpha,\beta,\gamma,\delta)^t\bigr)_{ij}&=\begin{cases}[u^iv^j]\,
\frac{1}{1-\alpha(u+v)-\alpha^2(\beta^2-1)uv},&j\le n-2\\[2mm]
[u^i]\,\frac{\gamma}{1-\alpha(1+\beta\delta)u},&j=n-1,\end{cases}\\[2mm]
\label{Ldetref}\det K(\alpha,\beta,\gamma,\delta)&=(\alpha\beta)^{(n-1)(n-2)/2}\,\gamma\,\delta^{n-1},\\
\label{LLKK}\makebox[0pt][c]{\hspace*{30mm}$\Bigl(L\bigl(\alpha_1,\beta_1\bigr)^{-1}\,L\bigl(\alpha_2,\beta_2\bigr)\,
\Bigl(K\bigl(\alpha_1,\beta_1,\gamma_1,\frac{\alpha_1\beta_1(\gamma_2-\gamma_1)}{\gamma_2(\alpha_2-\alpha_1)}\bigr)^{-1}\,
K\bigl(\alpha_2,\beta_2,\gamma_2,\frac{\alpha_2\beta_2(\gamma_2-\gamma_1)}{\gamma_1(\alpha_2-\alpha_1)}\bigr)\Bigr)^t\Bigr)_{ij}=$}\\
\notag\makebox[0pt][c]{\hspace*{65mm}$\begin{cases}\sum_{k=0}^{\min(i,j)}{i\choose k}{j\choose k}
\bigl(\frac{\alpha_2-\alpha_1}{\alpha_1\,\beta_1}\bigr)^{i+j}\,
\bigl(\frac{\alpha_2\,\beta_2}{\alpha_2-\alpha_1}\bigr)^{2k},&j\le n-2\\
\sum_{k=0}^i\sum_{l=0}^k{i\choose k}{n-l-2\choose k-l}\,
\bigl(\frac{\alpha_2-\alpha_1}{\alpha_1\,\beta_1}\bigr)^{i+n-1}
\bigl(\frac{\alpha_2\,\beta_2}{\alpha_2-\alpha_1}\bigr)^{2k}\,
\bigl(\frac{\gamma_2}{\gamma_1}\bigr)^{l+1},&j=n-1,\end{cases}$}\\
\notag\makebox[0pt][c]{\hspace*{135mm}for any $\alpha_1$, $\beta_1$, $\gamma_1$, $\alpha_2$, $\beta_2$, $\gamma_2$,}
\end{align}
where, as before, $[u^i]f(u)$ and $[u^iv^j]f(u,v)$ denote the coefficients of
$u^i$ and $u^iv^j$ in power series $f(u)$ and $f(u,v)$ respectively, and
$L(\alpha,\beta)$ is given by~\eqref{L}.  Property~\eqref{LLKK} can be obtained by first showing that
\begin{multline}\label{KK}\ts\Bigl(K\bigl(\alpha_1,\beta_1,\gamma_1,\frac{\alpha_1\beta_1(\gamma_2-\gamma_1)}{\gamma_2(\alpha_2-\alpha_1)}\bigr)^{-1}\,
K\bigl(\alpha_2,\beta_2,\gamma_2,\frac{\alpha_2\beta_2(\gamma_2-\gamma_1)}{\gamma_1(\alpha_2-\alpha_1)}\bigr)\Bigr)_{ij}=\\
\begin{cases}{i\choose j}\,
\bigl(\frac{\alpha_2-\alpha_1}{\alpha_1\,\beta_1}\bigr)^i\,
\bigl(\frac{\alpha_2\,\beta_2}{\alpha_2-\alpha_1}\bigr)^j,&i\le n-2\\[3mm]
\bigl(\frac{\alpha_2-\alpha_1}{\alpha_1\,\beta_1}\bigr)^{n-j-1}\,
\bigl(\frac{\alpha_2\,\beta_2}{\alpha_1\,\beta_1}\bigr)^j
\sum_{k=0}^j{n-k-2\choose j-k}\,\bigl(\frac{\gamma_2}{\gamma_1}\bigr)^{k+1},&i=n-1,\end{cases}\end{multline}
and using~\eqref{LL}.  Properties~\eqref{LK} and \eqref{KK} can be obtained using Cauchy-type residue integrals, similarly to the derivation
of~\eqref{LLTGF} and \eqref{LL} as indicated in \eqref{LLderiv}.

We now have
\begin{alignat*}{2}
&\ts\bigl(1+\omega_+\bigl(\frac{\tilde{a}\,b}{\rule{0ex}{1.65ex}a\,\tilde{b}}-1\bigr)\bigr)\,
Z_\ASM\bigl(n,\bigl(\frac{a}{b}\bigr)^2,\bigl(\frac{c}{b}\bigr)^2,\frac{\tilde{a}\,b}{\rule{0ex}{1.65ex}a\,\tilde{b}}\bigr)\\
&\ts=\frac{b}{\rule{0ex}{1.65ex}\tilde{b}}\,Z_\ASM\bigl(n,\bigl(\frac{a}{b}\bigr)^2,\bigl(\frac{c}{b}\bigr)^2,
\frac{\tilde{a}\,b}{\rule{0ex}{1.65ex}a\,\tilde{b}}\bigr)&\llap{[using \eqref{wr}, \eqref{nuqr} and \eqref{wrs}]}\\[2mm]
&\ts=\frac{1}{\rule{0ex}{1.65ex}b^{n(n-2)}\,\tilde{b}^n\,c^{n-1}\,\tilde{c}}\,Z_\SVDWBC(s,r,\ldots,r;r,\ldots,r)&\llap{[using \eqref{ZCZASMref}]}\\
&\ts=\frac{r^{n(n+1)}\,(\tilde{a}b)^n\,a^{n(n-1)}}{(r-s)^{n-1}\,c^n}\ds\det_{0\le i,j\le n-1}\!
\left(\begin{cases}[v^j]\left(\frac{1}{s(v+r)-q^2}-\frac{1}{s(v+r)-q^{-2}}\right),&i=0\\[2mm]
[u^{i-1}v^j]\left(\frac{1}{(u+r)(v+r)-q^2}-\frac{1}{(u+r)(v+r)-q^{-2}}\right),&i\ge1\end{cases}\right)\\
&&\llap{[using \eqref{Idet2}, \eqref{wr}, \eqref{DD}, \eqref{wrs} and \eqref{DDdetref}]}\\
&\ts=\frac{r^{n(n+1)}\,(\tilde{a}b)^n\,a^{n(n-1)}}{(s-r)^{n-1}\,c^n}\ds\det_{0\le i,j\le n-1}\!
\left(\begin{cases}[u^iv^j]\left(\frac{1}{(u+r)(v+r)-q^2}-\frac{1}{(u+r)(v+r)-q^{-2}}\right),&j\le n-2\\[2mm]
[u^i]\left(\frac{1}{s(u+r)-q^2}-\frac{1}{s(u+r)-q^{-2}}\right),&j=n-1\end{cases}\right)\\
&&\llap{[cycling the rows of and transposing the matrix]}\\
&\ts=\frac{r^{n(n+1)}\,(\tilde{a}b)^n\,a^{n(n-1)}}{(s-r)^{n-1}\,c^n}\,
\det\Bigl(\frac{1}{r^2-q^2}\,L\bigl(\frac{r}{q^2-r^2},\frac{q}{r}\bigr)\,
K\bigl(\frac{r}{q^2-r^2},\frac{q}{r},\frac{r^2-q^2}{rs-q^2},\frac{r-s}{q^{-1}rs-q}\bigr)^t\;-\\
&&\llap{$\frac{1}{r^2-q^{-2}}\,L\bigl(\frac{r}{q^{-2}-r^2},\frac{1}{qr}\bigr)\,
K\bigl(\frac{r}{q^{-2}-r^2},\frac{1}{qr},\frac{r^2-q^{-2}}{rs-q^{-2}},\frac{r-s}{qrs-q^{-1}}\bigr)^t\Bigr)$\ \ \text{[using \eqref{LK}]}}\\
&\ts=\frac{r^{n(n+1)}\,(\tilde{a}b)^n\,a^{n(n-1)}}{(qr^2-q^{-1})^{n(n-2)}\,(qrs-q^{-1})^n\,\,c^n}\,
\det\Bigl(\frac{1}{r^2-q^2}\,L\bigl(\frac{r}{q^{-2}-r^2},\frac{1}{qr}\bigr)^{-1}L\bigl(\frac{r}{q^2-r^2},\frac{q}{r}\bigr)\times\\
&&\llap{$\Bigl(K\bigl(\frac{r}{q^{-2}-r^2},\frac{1}{qr},\frac{r^2-q^{-2}}{rs-q^{-2}},\frac{r-s}{qrs-q^{-1}}\bigr)^{-1}
K\bigl(\frac{r}{q^2-r^2},\frac{q}{r},\frac{r^2-q^2}{rs-q^2},\frac{r-s}{q^{-1}rs-q}\bigr)\Bigr)^t\,-\,\frac{1}{r^2-q^{-2}}\,I\Bigr)$\quad}\\
&&\llap{[using \eqref{Ldet} and \eqref{Ldetref}]}\\
&=\det_{0\le i,j\le n-1}\ts\!\left((1-\omega_+)\delta_{ij}+\omega_+\begin{cases}\sum_{k=0}^{\min(i,j)}{i\choose k}{j\choose k}\,
\bigl(\frac{c}{b}\bigr)^{i+j}\bigl(\frac{a}{c}\bigr)^{2k},&j\le n-2\\[3mm]
\sum_{k=0}^i\sum_{l=0}^k{i\choose k}{n-l-2\choose k-l}\,
\bigl(\frac{c}{b}\bigr)^{i+n-1}\bigl(\frac{a}{c}\bigr)^{2k}
\bigl(\frac{\tilde{a}\,b}{\rule{0ex}{1.65ex}a\,\tilde{b}}\bigr)^{l+1},&j=n-1\end{cases}\right)\\
&&\llap{[using \eqref{wr}, \eqref{nuqr}, \eqref{wrs} and \eqref{LLKK}]}\\
&=\det_{0\le i,j\le n-1}\ts\!\left((1-\omega_+)\delta_{ij}+\omega_+\begin{cases}\sum_{k=0}^{\min(i,j)}{i\choose k}{j\choose k}\,
\bigl(\frac{a}{b}\bigr)^{2k}\bigl(\frac{c}{b}\bigr)^{2(i-k)},&j\le n-2\\[3mm]
\sum_{k=0}^i\sum_{l=0}^k{i\choose k}{n-l-2\choose k-l}\,
\bigl(\frac{a}{b}\bigr)^{2k}\bigl(\frac{c}{b}\bigr)^{2(i-k)}
\bigl(\frac{\tilde{a}\,b}{\rule{0ex}{1.65ex}a\,\tilde{b}}\bigr)^{l+1},&j=n-1\end{cases}\right)\\
&&\llap{[multiplying the matrix on the left by $V$ and right by $V^{-1}$, where $V_{ij}=\bigl(\frac{c}{b}\bigr)^i\delta_{ij}$]}.
\end{alignat*}
Applying the transformation $r\rightarrow r^{-1}$,~$s\rightarrow s^{-1}$ and~$q\rightarrow q^{-1}$
to the equality between the first and last expressions in the previous sequence,
it follows that this equality remains valid if~$\omega_+$ is replaced by~$\omega_-$.
Finally, identifying $x=\bigl(\frac{a}{b}\bigr)^2$, $y=\bigl(\frac{c}{b}\bigr)^2$ and
$z=\frac{\tilde{a}\,b}{\rule{0ex}{1.65ex}a\,\tilde{b}}$ gives the result of Proposition~\ref{ASMdetref},
using \eqref{nuxy}, \eqref{nuabc} and \eqref{MASMref}.

\subsection{\protect Proof of Proposition \ref{DPPdetref}}\label{DPPdetderivref}
In this subsection, the DPP determinant formula of Proposition \ref{DPPdetref} will be obtained using
the same approach as in Section~\ref{DPPdetderiv},
i.e., the Lindstr\"{o}m--Gessel--Viennot theorem~\eqref{LGV} and the bijection between DPPs and
certain sets of nonintersecting lattice paths will be used.

\begin{figure}[h]
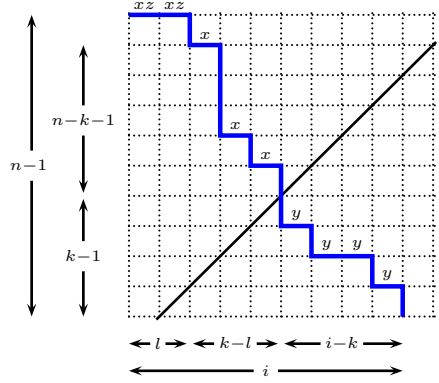
\centering
\psset{unit=4mm}
\pspicture(-3.5,-1.1)(10,10)
\psgrid[subgriddiv=1,griddots=6,gridlabels=0](0,0)(10,10)
\psline[linewidth=1pt](0.9,-0.1)(10.1,9.1)
\rput(-1.5,2){$\sss k-1$}
\psline[linewidth=1pt]{<-}(0,-0.9)(0.7,-0.9)\rput(0.95,-0.9){$\sss l$}\psline[linewidth=1pt]{->}(1.2,-0.9)(1.9,-0.9)
\psline[linewidth=1pt]{<-}(2.1,-0.9)(2.8,-0.9)\rput(3.5,-0.9){$\sss k-l$}\psline[linewidth=1pt]{->}(4.2,-0.9)(4.9,-0.9)
\psline[linewidth=1pt]{<-}(5.1,-0.9)(6.3,-0.9)\rput(7,-0.9){$\sss i-k$}\psline[linewidth=1pt]{->}(7.7,-0.9)(9,-0.9)
\psline[linewidth=1pt]{<-}(0,-1.8)(4.1,-1.8)\rput(4.5,-1.8){$\sss i$}\psline[linewidth=1pt]{->}(4.9,-1.8)(9,-1.8)
\psline[linewidth=1pt]{<-}(-1.5,0)(-1.5,1.5)\psline[linewidth=1pt]{->}(-1.5,2.5)(-1.5,3.9)
\psline[linewidth=1pt]{<-}(-1.5,4.1)(-1.5,6)\rput(-1.5,6.5){$\sss n-k-1$}\psline[linewidth=1pt]{->}(-1.5,7)(-1.5,9)
\psline[linewidth=1pt]{<-}(-3.2,0)(-3.2,4.5)\rput(-3.3,5){$\sss n-1$}\psline[linewidth=1pt]{->}(-3.2,5.5)(-3.2,10)
\psline[linewidth=1.7pt,linecolor=blue](0,10)(2,10)(2,9)(3,9)(3,6)(4,6)(4,5)(5,5)(5,3)(6,3)(6,2)(8,2)(8,1)(9,1)(9,0)
\rput[b](0.5,10.2){$\sss xz$}\rput[b](1.5,10.2){$\sss xz$}\rput[b](2.5,9.2){$\sss x$}\rput[b](3.5,6.2){$\sss x$}\rput[b](4.5,5.2){$\sss x$}
\rput[b](5.5,3.2){$\sss y$}\rput[b](6.5,2.2){$\sss y$}\rput[b](7.5,2.2){$\sss y$}\rput[b](8.5,1.2){$\sss y$}
\endpspicture
\caption{Derivation of~\eqref{DPPwpref} for $j=n-1$.}\label{DPPdetfigref}\end{figure}

Consider again the set $\NILP(n)$, as defined in Section~\ref{DPPdetderiv}.
It follows from the bijection between $\DPP(n)$ and $\NILP(n)$ that the statistic~$\rho$ for DPPs, as given by~\eqref{rhoD},
satisfies the following result.
\begin{lemma}\label{rhoDPP}
Let $D\in\DPP(n)$ correspond to $P\in\NILP(n)$.  Then
\[\rho(D)=\text{number of rightward steps in $P$ in the top row of the grid.}\]
\end{lemma}
Now assign a weight of~$xz$ to each horizontal edge in the top row of the grid, and assign the same weights as in Section~\ref{DPPdetderiv}
to all other edges.
Then the bijection between $\DPP(n)$ and~$\NILP(n)$,
together with~\eqref{ZDPP}, \eqref{NILP}, \eqref{detidentity} and Lemmas~\ref{numuDPP} and~\ref{rhoDPP}, gives
\begin{equation}\label{DPPLGVref}
Z_\DPP(n,x,y,z)=
\det_{0\le i,j\le n-1}\Biggl(-\delta_{i,j+1}+\sum_{p\in\paths((0,j),(i,0))}W(p)\Biggr).\end{equation}
For $j\le n-2$, each path of
$\paths((0,j),(i,0))$ remains below the top row of the grid, and
so~\eqref{DPPwp} holds. For
$j=n-1$, the paths of $\paths((0,n-1),(i,0))$ can be obtained by
joining the unique path of $\paths((0,n-1),(l,n-1))$, which has
weight $(xz)^l$, the unique path of $\paths((l,n-1),(l,n-2))$, which
has weight~$1$, any of the $n-l-1\choose k-l$ paths of
$\paths((l,n-2),(k,k-1))$, each of which has weight~$x^{k-l}$, and
any of the $i-1\choose i-k$ paths of $\paths((k,k-1),(i,0))$, each
of which has weight~$y^{i-k}$, for any $0\le k\le i$ and $0\le l\le k$.
This is shown diagrammatically in Figure~\ref{DPPdetfigref}.
It therefore follows that
\begin{equation}\label{DPPwpref}
\sum_{p\in\paths((0,j),(i,0))}W(p)=
\begin{cases}
\ds\sum_{k=0}^{\min(i,j+1)}{i-1\choose i-k}{j+1\choose k}x^k\,y^{i-k},&j\le n-2\\
\ds\sum_{k=0}^i\,\sum_{l=0}^k{i-1\choose i-k}{n-l-1\choose k-l}x^k\,y^{i-k}\,z^l,&j=n-1.
\end{cases}
\end{equation}

Finally, combining~\eqref{DPPLGVref} and~\eqref{DPPwpref},
and multiplying the last column of the matrix by~$1+\omega(z-1)$, gives
the result of Proposition~\ref{DPPdetref}.

The alternative sets of nonintersecting lattice paths described in
Section~\ref{DPPdetderiv} and the TB configurations used by Lalonde~\cite{Lal-ASM,Lal-DPP}
lead to the alternative formulae
$Z_\DPP(n,x,y,z)=\det M'_\DPP(n,x,y,z)=\det M''_\DPP(n,x,y,z)$, where
the entries of the matrices here are given, for
$0\le i,j\le n-1$, by
$M'_\DPP(n,x,y,z)_{ij}=M'_\DPP(n,x,y,1)_{ij}$ and $M''_\DPP(n,x,y,z)_{ij}=
M''_\DPP(n,x,y,1)_{ij}$ if $j\le n-2$ (with $M'_\DPP(n,x,y,1)$ and $M''_\DPP(n,x,y,1)$
given by~\eqref{Md} and~\eqref{Mdd}), and
\begin{align}\label{Mdref}M'_\DPP(n,x,y,z)_{i,n-1}&=\delta_{i,n-1}+
\sum_{k=0}^{i-1}\,\sum_{l=0}^k\,\sum_{m=0}^l{n-m-2\choose l-m}{k\choose l}x^{l+1}\,y^{k-l}\,z^{m+1},\\
\label{Mddref}M''_\DPP(n,x,y,z)_{i,n-1}&=
\sum_{k=0}^i{n-k-1\choose i-k}x^i\,z^k.\end{align}
Note that $(I-S)M'_\DPP(n,x,y,z)=\bar{M}(n,x,y,z)(I-S^t)$ and that $B(y)M''_\DPP(n,x,y,z)=\bar{M}(n,x,y,z)$,
where~$\bar{M}(n,x,y,z)$ (as will be given explicitly in~\eqref{Mbar}) is the matrix obtained from $M_\DPP(n,x,y,z)$ by omitting the
factor $1+\omega(z-1)$ from the last column, $S$ is given by~\eqref{S} and $B(y)_{ij}={i-1\choose i-j}\,y^{i-j}$.
The matrix~$M'_\DPP(n,x,y,z)$ is given, using slightly different notation, in Mills, Robbins and Rumsey~\cite[p.~346]{MRR-ASM}, and
the result $Z_\DPP(n,x,y,z)\!=\det M'_\DPP(n,x,y,z)$ is stated there, without proof.

\subsection{\protect Proof of Proposition \ref{ASMdetDPPdetref}}\label{ASMdetDPPdetproofref}
In this subsection, the equality of the ASM and DPP determinants in Proposition~\ref{ASMdetDPPdetref} will be proved using
the same approach as in Section~\ref{ASMdetDPPdetproof}, i.e.,
a certain relation between the generating functions associated with the entries of the
respective matrices will be obtained, and this will then imply a relation between the matrices
and the equality of their determinants.

The entries of the last columns of $M_\ASM(n,x,y,z)$ and $M_\DPP(n,x,y,z)$ satisfy
\begin{align}\label{MASMlast}M_\ASM(n,x,y,z)_{i,n-1}&=
(1-\omega)\delta_{i,n-1}+[u^i]\,\frac{\omega z}{1-(xz-x+y)u}\left(1+\frac{xu}{1-yu}\right)^n,\\[2mm]
\label{MDPPlast}M_\DPP(n,x,y,z)_{i,n-1}&=[u^i]\,\frac{(1+\omega(z-1))(1-yu)}{1-(xz-x+y)u}\left(1+\frac{xu}{1-yu}\right)^n.\end{align}
Now define the generating functions
\begin{align}\ds f_\ASM(n,x,y,z;u,v)=&\frac{1-\omega}{1-uv}+\frac{\omega}{1-yu-v-(x-y)uv}\\[2mm]
\notag&&\llap{$\ds+\,\frac{\omega(z-1)}{1-(xz-x+y)u}\left(1+\frac{xu}{1-yu}\right)^nv^{n-1}
\left(1+\frac{(\omega x-\omega-x+\omega yu)v}{y(\omega+xu-\omega y u)}\right)$},\\[2mm]
\ds f_\DPP(n,x,y,z;u,v)=&\frac{1-u}{(1-v)(1-uv)}+\frac{xu}{(1-v)(1-yu-v-(x-y)uv)}\\[2mm]
\notag&\ds+\,\frac{(z-1)(\omega+xu-\omega yu)}{1-(xz-x+y)u}\left(1+\frac{xu}{1-yu}\right)^nv^{n-1}.\end{align}
Then, using \eqref{gASM}, \eqref{gDPP}, \eqref{MASMlast} and \eqref{MDPPlast}, it follows that
\begin{align}\label{fASMref}M_\ASM(n,x,y,z)_{ij}&=[u^iv^j]f_\ASM(n,x,y,z;u,v),\\[2mm]
M_\DPP(n,x,y,z)_{ij}&=[u^iv^j]f_\DPP(n,x,y,z;u,v).\end{align}
Using the equation~\eqref{nuxy} satisfied by $\omega$, it is seen that
\begin{equation}\label{fASMMfDPPref}(1+(x-\omega y-1)u)f_\ASM(n,x,y,z;u,v)=(1+(\omega-1)v)f_\DPP(n,x,y,z;u,v).\end{equation}
(In particular, the difference between the two sides of~\eqref{fASMMfDPPref} contains
the LHS of~\eqref{nuxy} as a factor, and therefore vanishes.)
Note that if the final term $1+\frac{(\omega x-\omega-x+\omega yu)v}{y(\omega+xu-\omega y u)}$ were omitted from
$f_\ASM(n,x,y,z;u,v)$, then~\eqref{fASMref} would still hold, but \eqref{fASMMfDPPref} would not.

It now follows that
\begin{equation}\label{MASMMDPPref}
(I+(x-\omega y-1)S)\,M_\ASM(n,x,y,z)=M_\DPP(n,x,y,z)\,(I+(\omega-1)S^t),\end{equation}
where $S$ is again given by~\eqref{S}, and this implies the result of Proposition~\ref{ASMdetDPPdetref}.

\section{Discussion}\label{discuss}
The proof of the Mills, Robbins and Rumsey ASM-DPP conjecture (i.e., Theorem~\ref{MRRC1}),
which has been the primary focus of this paper, has been presented in the preceding two sections.
In this final section, we discuss some consequences for the computation of the associated ASM and DPP generating functions,
describe certain symmetry operations for ASMs and DPPs, outline some further work which generalizes the
work reported in this paper,
and review some remaining unresolved matters regarding the relationship between ASMs and DPPs.

\subsection{Determinant formulae for weighted enumeration of ASM and DPPs}
This paper has mainly been devoted to the proof of Theorem~\ref{MRRC1},
which asserts the equality of sizes of subsets of $\ASM(n)$ and $\DPP(n)$ associated with arbitrary values of certain statistics,
or equivalently the equality of the generating functions which give weighted enumerations of the elements of $\ASM(n)$ or~$\DPP(n)$
using arbitrary weights associated with these statistics.  However, as a byproduct of this proof, determinant formulae
have been obtained which enable the efficient computation of these generating functions for given $n$, and thereby
the evaluation of the sizes of all of the associated subsets of $\ASM(n)$ and $\DPP(n)$.
In particular, it follows from Propositions~\ref{ASMdetref}--\ref{ASMdetDPPdetref} that
the generating functions \eqref{ZASM} and \eqref{ZDPP} satisfy
\begin{equation}\label{wtenum}
Z_\ASM(n,x,y,z)=Z_\DPP(n,x,y,z)=\det\bar{M}(n,x,y,z),
\end{equation}
where
\begin{equation}\label{Mbar}
\bar{M}(n,x,y,z)_{ij}=\!-\delta_{i,j+1}+
\begin{cases}\ds\sum_{k=0}^{\min(i,j+1)}{i-1\choose i-k}{j+1\choose k}x^k\,y^{i-k},&j\le n-2\\
\ds\sum_{k=0}^i\,\sum_{l=0}^k{i-1\choose i-k}{n-l-1\choose k-l}x^k\,y^{i-k}\,z^l,&j=n-1,\end{cases}
\end{equation}
with $0\le i,j\le n-1$.

Note that the matrices $M'_\DPP(n,x,y,z)$ or
$M''_\DPP(n,x,y,z)$, as given by~\eqref{Md}--\eqref{Mdd} and~\eqref{Mdref}--\eqref{Mddref}, could be used here instead
of $\bar{M}(n,x,y,z)$,
with the formulae for $Z_\DPP(n,x,y,z)$ as the determinant of these matrices having been
obtained by Mills, Robbins and Rumsey~\cite[p.~346]{MRR-ASM} and Lalonde~\cite[Thm.~3.1]{Lal-ASM} respectively.

It can be seen that $Z_\ASM(n,x,y,z)=Z_\DPP(n,x,y,z)=\frac{1}{1+\omega(z-1)}\det M_\ASM(n,x,y,z)$ also holds,
where $M_\ASM(n,x,y,z)$ is given by~\eqref{MASMref} and~$\omega$ is a solution of~\eqref{nuxy}.
However, this form seems computationally less useful, due to the presence in the prefactor and the matrix of the
term~$\omega$, which does not have polynomial dependence on $x$, $y$ and~$z$.
A formula for $Z_\ASM(n,x,y,1)$ in terms of the determinant of a matrix closely related to $M_\ASM(n,x,y,1)$ was also obtained by Colomo and
Pronko in~\cite[Eqs.~(23)--(24)]{CP-6v0} and~\cite[Eqs.~(4.3)--(4.7)]{CP-6v},
but it seems that a general determinant formula
for $Z_\ASM(n,x,y,z)$ of the form of~\eqref{wtenum}--\eqref{Mbar}
(i.e., in which all entries of the matrix are polynomials in $x$, $y$ and $z$) has not previously appeared in the literature.
In fact, the equation $Z_\ASM(n,1,y,z)=\det M'_\DPP(n,1,y,z)$, with $M'_\DPP(n,1,y,z)$ given by~\eqref{Md}--\eqref{Mdd},
confirms the validity of a determinant formula conjectured by Robbins~\cite[Conj.~3.1]{Rob-ASM}.
(See also de Gier~\cite[Sec.~2]{dG-review3}.)

Note also that combining~\eqref{ZCZASM},~\eqref{wtenum} and~\eqref{Mbar} gives a determinant formula for the
partition function of the homogeneous six-vertex model with DWBC directly
in terms of the weights~$a$,~$b$ and $c$, rather than in terms of the spectral parameters,
as is the case with previously-known determinant formulae for this partition function.
(Specifically, this formula for an $n\times n$ grid is $c^n\det_{0\le i,j\le n-1}\bigl(-b^{2i}\delta_{i,j+1}+
\sum_{k=0}^{\min(i,j+1)}{i-1\choose i-k}{j+1\choose k}a^{2k}c^{2(i-k)}\bigr)$.)

Finally, note that it seems that the techniques of this paper could be used to obtain
similar such determinant formulae for the partition functions of cases of the homogeneous six-vertex model
with certain other boundary conditions, for example some of the cases studied by Kuperberg in~\cite{Kup-symASM}.
Indeed, the study of these cases in~\cite{Kup-symASM} led to enumerative results for various classes of ASMs,
and, using the Lindstr\"{o}m--Gessel--Viennot theorem,
determinant formulae of the form of~\eqref{wtenum} can often be interpreted in terms of the enumeration of particular
types of plane partition.  Accordingly, if such determinant formulae are obtained for these classes of ASMs, then
this may reveal new enumerative connections between these ASMs and certain plane partitions.

\subsection{Symmetry operations}\label{symm}
We now discuss certain symmetry operations for ASMs and DPPs,
and state a result, analogous to Theorem~\ref{MRRC1}, for the equality
of sizes of subsets of $\ASM(n)$ and $\DPP(n)$ which are invariant under these operations.

Define an operation $\ast$ on $\ASM(n)$ and on $\DPP(n)$ as follows. For $A\in\ASM(n)$, let~$A^\ast$ be the
ASM obtained by reflecting~$A$ in a vertical line passing through the center of the matrix, i.e., $A^*_{ij}=A_{i,n+1-j}$.
For DPPs, consider the bijection, as obtained by Krattenthaler~\cite{Kratt-DPP} and shown for the example of~\eqref{DPPEx} in Figure~\ref{hex},
between $\DPP(n)$ and the set of cyclically symmetric rhombus tilings of a hexagon
with alternating sides of length $n-1$ and~$n+1$ and a central equilateral triangular hole of side length~2.
If $D\in\DPP(n)$ corresponds to the rhombus tiling~$R$,
then~$D^*$ is defined to be the DPP which corresponds to the reflection of~$R$ in any of the three lines
which bisect the central triangular hole.

For example, for the DPP~$D$ of~\eqref{DPPEx}, $D^\ast=\ba{c@{\;}c@{\;}c@{\;}c@{\;}c}6&6&5&5&3\\[-1mm]
&4&2&2\ea$ (taking $n=6$).  As further examples, for~$\ASM(3)$ and~$\DPP(3)$, as given in~\eqref{ASM3} and~\eqref{DPP3},
the first and second, third and fourth, and fifth and sixth elements are
related by $\ast$, and the seventh elements are invariant under~$\ast$.

It can immediately be seen that~$\ast$ is an involution on $\ASM(n)$ and on $\DPP(n)$.
Mills, Robbins and Rumsey~\cite[Sec.~3]{MRR-ASM} showed that~$\ast$ is also a unique antiautomorphism on
$\DPP(n)$, with respect to a certain partial order.

The operation~$\ast$ for DPPs was introduced by Mills, Robbins and Rumsey~\cite[Sec.~3]{MRR-ASM},
and first defined directly in terms of the parts of a DPP~\cite[p.~351]{MRR-ASM}.  (See also Bressoud~\cite[pp.~194--195]{Bressoud}.)
The alternative description of~$\ast$ for DPPs in terms of rhombus tilings was first obtained by Krattenthaler~\cite[p.~1143]{Kratt-DPP}.
A description of~$\ast$ for DPPs in terms of the TB configurations of Lalonde is given in~\cite{Lal-DPP},
and a closely-related description can be given in terms of the set~$\NILP(n)$ of sets of nonintersecting lattice paths,
as defined in Section~\ref{DPPdetderiv}.

It was shown by Mills, Robbins and Rumsey~\cite[p.~352]{MRR-ASM}
(see also Krattenthaler~\cite{Kratt-DPP} and Lalonde~\cite{Lal-DPP} for the DPP case) that
the ASM and DPP statistics~\eqref{nuA}--\eqref{rhoD} behave under~$\ast$ according to
\begin{equation}\label{symstat}\ts\nu(X^\ast)=\frac{n(n-1)}{2}-\nu(X)-\mu(X),\qquad\mu(X^\ast)=\mu(X),\qquad
\rho(X^\ast)=n-1-\rho(X),\end{equation}
for any $X\in\ASM(n)$ or $X\in\DPP(n)$.
(It can also be shown that, for any $D\in\DPP(n)$, the sum of the number of rows in $D$ and
number of rows in $D^\ast$ is $n-1$, and the sum of all parts in both $D$ and $D^\ast$ is $2{n+1\choose 3}$.)

The following result was conjectured by Mills, Robbins and Rumsey~\cite[Conj.~3S]{MRR-ASM}
and proved (in a form involving slightly different terminology) by de Gier, Pyatov and Zinn-Justin~\cite[Prop.~3, first equation]{artic44}.
\begin{theorem*}[de Gier, Pyatov, Zinn-Justin]
The sizes of
$\{A\in\ASM(n)\mid A^\ast=A,\;\mu(A)=m\}$
and $\{D\in\DPP(n)\mid D^\ast=D,\;\mu(D)=m\}$
are equal for any $n$ and $m$.
\end{theorem*}
It follows from~\eqref{symstat} and other basic properties of the operations~$\ast$
that if the sets in the previous theorem are nonempty, then $n$ is odd, $m/2-(n-1)/4$ is a nonnegative integer,
and any element $X$ of either set satisfies
$\nu(X)=n(n-1)/4-m/2$ and $\rho(X)=(n-1)/2$.  As examples,
\begin{align}&\{A\in\ASM(5)\mid A^\ast=A\}=\\
\notag&\qquad\qquad\left\{\begin{pmatrix}0&0&1&0&0\\0&1&-1&1&0\\0&0&1&0&0\\1&0&-1&0&1\\0&0&1&0&0\end{pmatrix}\!,
\begin{pmatrix}0&0&1&0&0\\1&0&-1&0&1\\0&0&1&0&0\\0&1&-1&1&0\\0&0&1&0&0\end{pmatrix}\!,
\begin{pmatrix}0&0&1&0&0\\0&1&-1&1&0\\1&-1&1&-1&1\\0&1&-1&1&0\\0&0&1&0&0\end{pmatrix}\right\},\\[3mm]
&\{D\in\DPP(5)\mid D^\ast=D\}=\\[-2mm]
\notag&\qquad\qquad\left\{\ba{c@{\;}c@{\;}c@{\;}c}5&5&3&2\\[-0.8mm]&4&1\ea,\,
\ba{c@{\;}c@{\;}c@{\;}c}5&5&4&2\\[-0.8mm]&3&1\ea,\,
\ba{c@{\;}c@{\;}c@{\;}c}5&5&2&2\\[-0.8mm]&4&1&1\ea
\right\}\cong
\left\{\vcenter{\hbox{\includegraphics[width=2cm]{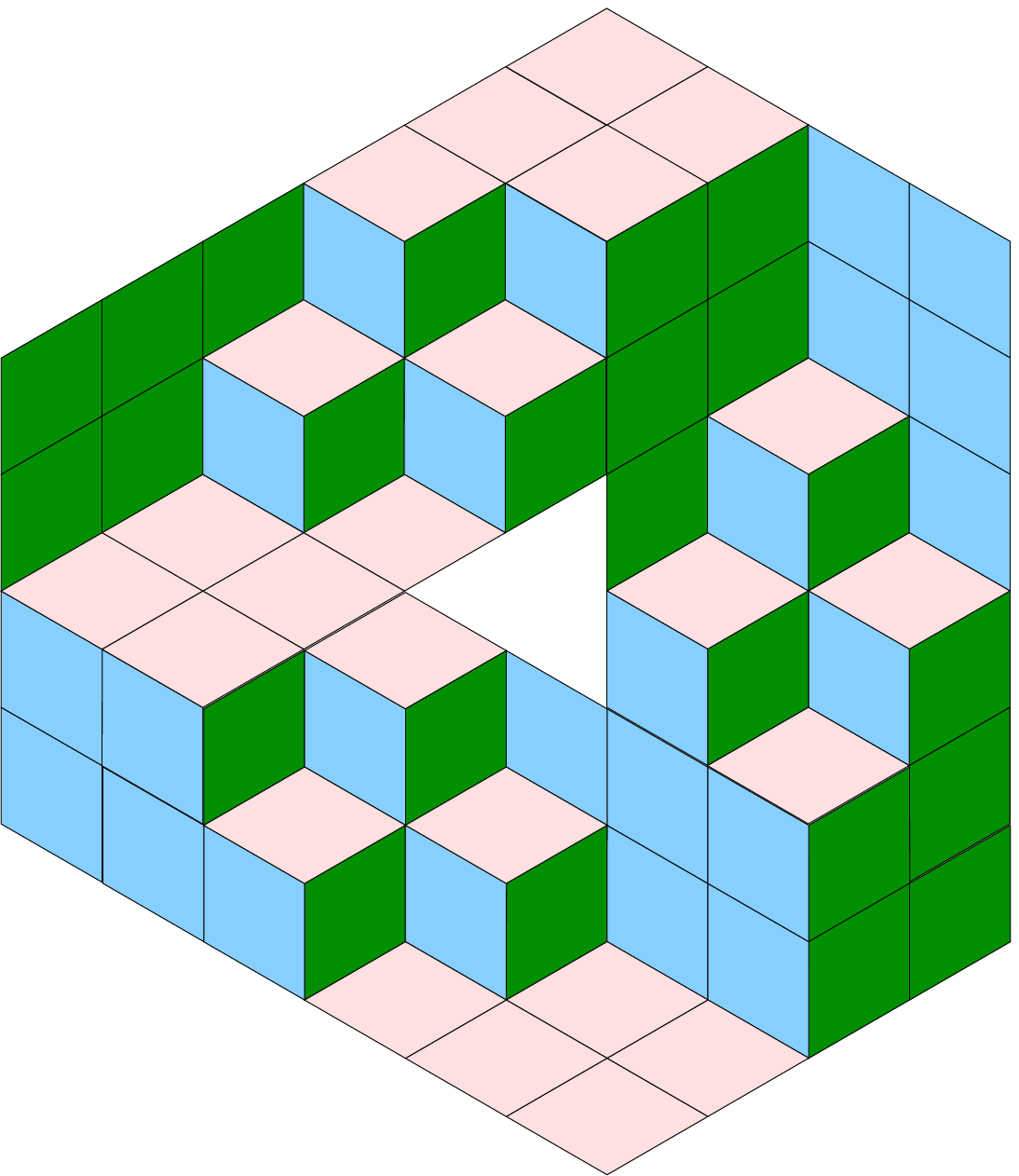}}},
\vcenter{\hbox{\includegraphics[width=2cm]{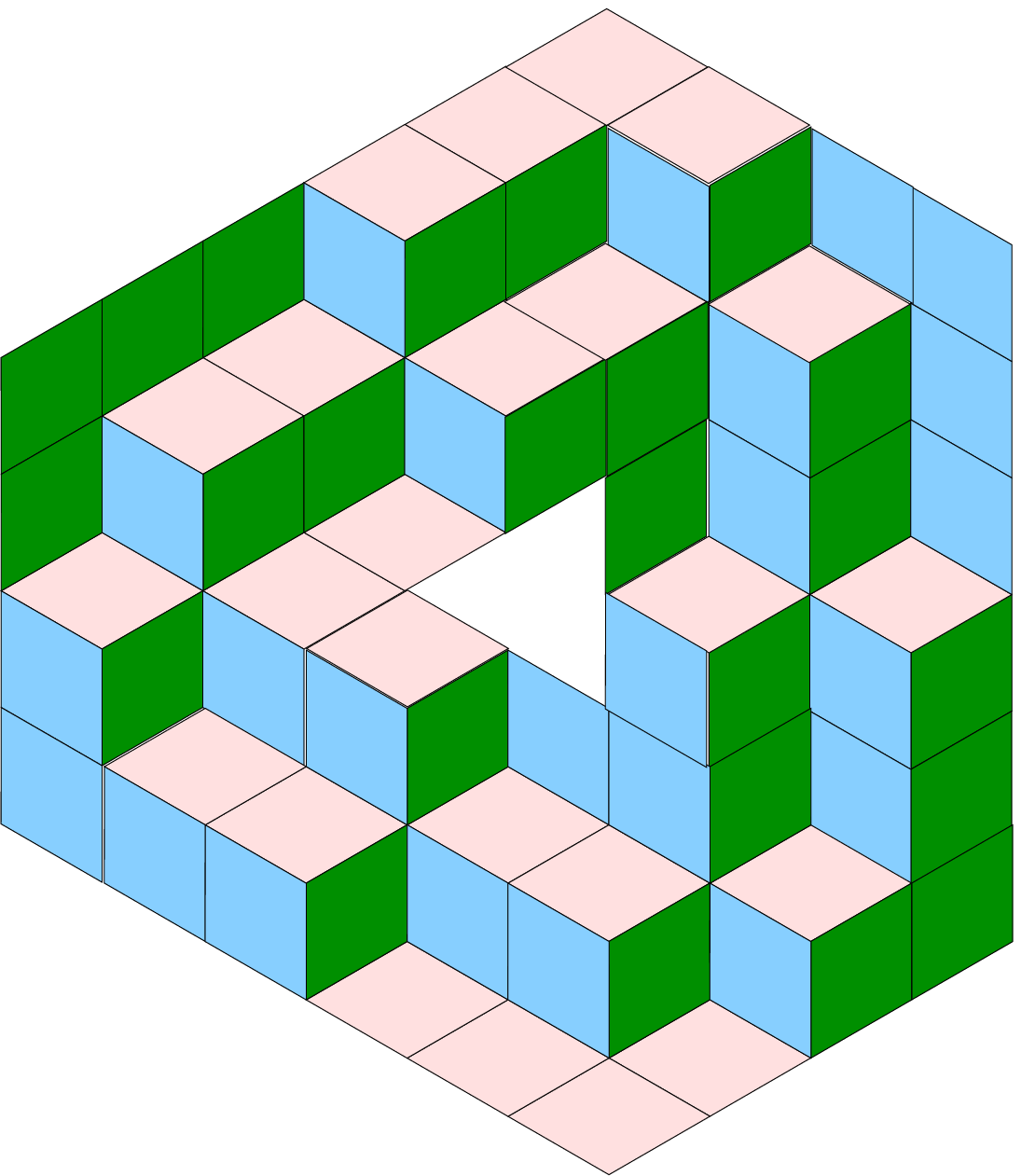}}},
\vcenter{\hbox{\includegraphics[width=2cm]{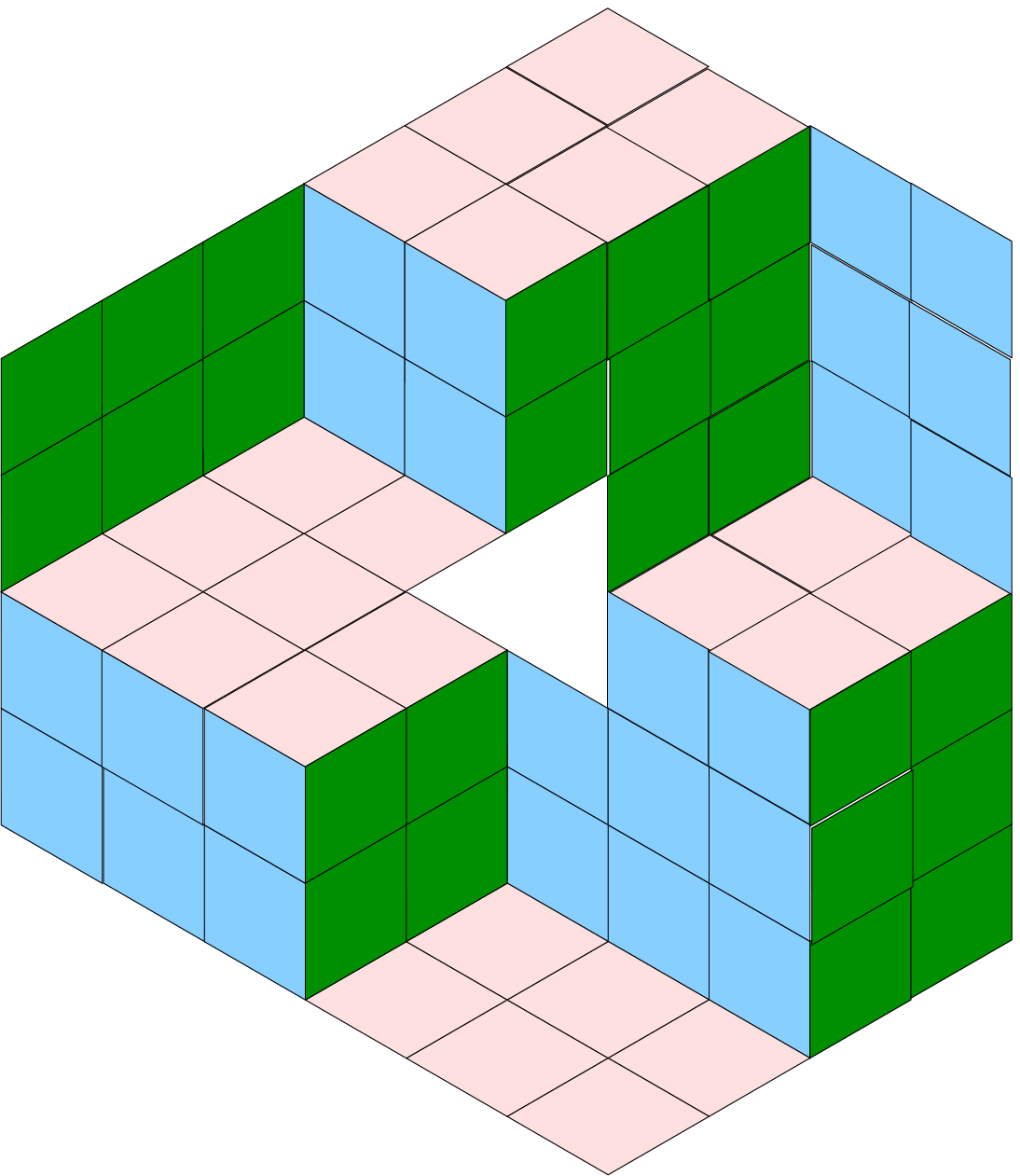}}}
\right\}.\end{align}

A bijective proof of the $m=0$ case of the previous theorem can be obtained straightforwardly.
Also, for the case of the theorem with $m$ summed over,
it is known that
\begin{equation}\label{SASMDPP}|\{A\in\ASM(2n+1)\mid A^\ast=A\}|=|\{D\in\DPP(2n+1)\mid
D^\ast=D\}|=\,\ds\prod_{i=1}^{n}\frac{(6i-2)!}{(2n+2i)!}.\end{equation}
The formula for ASMs in~\eqref{SASMDPP} was conjectured by Robbins~\cite{Rob-story,Rob-ASM} and
first proved by Kuperberg~\cite[Thm.~2, second equation]{Kup-ASM},
while the formula for DPPs can be obtained using the correspondence between DPPs and
sets of nonintersecting lattice paths, together with the Lindstr\"{o}m--Gessel--Viennot theorem, to obtain a determinant
which is then evaluated using a result of Mills, Robbins and Rumsey~\cite[\protect Thm.~7, 2nd part, $\mu=1$]{MRR-symPP}.
Certain aspects of $\ast$-invariant ASMs and DPPs are also considered by Di
Francesco and Zinn-Justin~\cite{DiFran,artic43}.

It has thus been seen in this subsection that the symmetry operations~$\ast$ on ASMs  and~$\ast$ on DPPs seem
to correspond in the sense that they
share the same behavior with regard to the statistics~$\nu$,~$\mu$ and~$\rho$
of~\eqref{nuA}--\eqref{rhoD}, and that the sizes of sets of $\ast$-invariant ASMs and DPPs are equal.
Accordingly, if a bijective proof of Theorem~\ref{MRRC1} is found, then the bijection would
preferably satisfy the property that for any ASM~$A$ and DPP~$D$ which correspond,~$A^\ast$ and~$D^\ast$ also correspond.

\subsection{Further work and remaining questions}
In this paper we have proved a result, Theorem~\ref{MRRC1}, which confirms the existence of a close relationship between
ASMs and DPPs.

In future work we plan to present the proof
of a certain doubly-refined generalization of Theorem~\ref{MRRC1} involving a further statistic for ASMs and DPPs
(where this statistic for an ASM is simply the number of $0$'s to the right of the $1$ in the last row),
to define further statistics and give associated results for $\ast$-invariant ASMs and DPPs (where~$\ast$ is the operation
defined in Section~\ref{symm}), and to discuss some progress made towards constructing a bijection between ASMs and DPPs.

However, despite the implications of our current and forthcoming results, there remain various questions about
the relationship between ASMs and DPPs which are likely only to be resolved by the discovery of a full bijection
between these two types of objects.
For example, there are certain features which, based on current knowledge, are displayed by only one type
of object, with no corresponding feature known for the other.  (Of course, even if a bijection is found,
it is likely that many of these features will still have a simple or natural description in terms of only one type of
object.)

A feature displayed by ASMs is the existence of eight natural symmetry operations,
corresponding to those of the dihedral group acting on a square matrix, i.e., the identity,
reflections in central vertical, horizontal, diagonal and antidiagonal lines, and rotations about the center by~$\pi/2$,~$\pi$ and~$3\pi/2$.
The operation $\ast$ for DPPs, as discussed in Section~\ref{symm}, seems to correspond to vertical reflection of ASMs, but
as was pointed out by Mills, Robbins and Rumsey~\cite[p.~353]{MRR-ASM} and is still the case, there are no
known operations for DPPs which seem to correspond to the last six of the ASM operations.

Another feature of ASMs is the existence of simple bijections with several other combinatorial objects, which themselves
display natural properties for which no analogous properties are known for DPPs.
A review of such bijections is given by Propp~\cite{Propp-ASM}.
For example, as described in~\cite{Propp-ASM}, the set $\ASM(n)$ is in bijection with the sets,
all associated with a certain order $n$, of monotone triangles,
height-function matrices and fully packed loop configurations.  The bijection with
monotone triangles enables a further statistic to be defined naturally for ASMs, where this statistic enters
a certain refined enumerative relationship between $\ASM(n)$ and the set of totally
symmetric self-complementary plane partitions in a $2n\times2n\times2n$ box
(see Mills, Robbins and Rumsey~\cite[Conj.~7]{MRR-TSSCPP} and Zeilberger~\cite[Lemma~1]{Zeil-ASM}),
the bijection with monotone triangles or
with height-function matrices enables $\ASM(n)$ to be regarded as a distributive lattice,
where this lattice is associated with various algebraic results (see for example
Lascoux~\cite[pp.~3--4]{Las-CY}, Lascoux and Sch\"{u}tzenberger~\cite{LS-Cox}, Propp~\cite[pp.~46--47]{Propp-ASM}
and Striker~\cite{Stri-Poset}),
and the bijection with fully packed loop configurations enables
ASMs to be classified according to certain link patterns, where these patterns play a primary role in the Razumov--Stroganov (ex-)conjecture
(see for example Cantini and Sportiello~\cite{CS-RS} and Propp~\cite[Sec.~7]{Propp-ASM}).

On the other hand, examples of some features which are exhibited by DPPs, but for which no counterparts are known for ASMs,
are the existence of a certain lattice structure, as described by Mills, Robbins and Rumsey~\cite[p.~347]{MRR-ASM}
(where this is a different lattice structure from the one known for ASMs),
and the existence of two natural statistics, given by the number of rows of a DPP and the sum of parts of a DPP.
Furthermore, it can be shown, using the derivations of Sections~\ref{DPPdetderiv} and~\ref{DPPdetderivref}, that
$\sum_{D\in\DPP(n)}w^{\phi(D)+1}\,x^{\nu(D)}\,y^{\mu(D)}\,z^{\rho(D)}=\det\bar{M}(n,w,x,y,z)$,
where~$\phi(D)$ is the number of rows of~$D$ and $\bar{M}(n,w,x,y,z)$ is the matrix obtained from $\bar{M}(n,x,y,z)$ by multiplying the second
term on the RHS of~\eqref{Mbar} by $w$, and it was shown by Mills, Robbins and Rumsey~\cite{MRR-Mac} that
$\sum_{D\in\DPP(n)}q^{|D|}=\prod_{i=0}^{n-1}\frac{[3i+1]_q!}{[n+i]_q!}$,
where~$|D|$ is the sum of parts of $D$ and $[n]_q!$ is defined after~\eqref{Z0}.

Finally, and intriguingly, it has been shown (nonbijectively) that certain enumerations of DPPs involving the sum of parts
are related to enumerations of ASMs invariant under certain symmetry operations.
More specifically, it follows from a result of Stanton, as given by Stephens-Davidowitz and
Cloninger~\cite[Thm.~2.2]{SDC-ASM}, that
\begin{align}&|\{D\in\DPP(n)\mid\,|D|\text{ even}\}|=|\{D\in\DPP(n)\mid\,|D|\text{ odd}\}|\;+\\
\notag&&\llap{$|\{A\in\ASM(n)\mid A \text{ invariant under rotation by }\pi\}|$},\\
&|\{D\in\DPP(n)\mid|D|\equiv0\bmod4\}|=|\{D\in\DPP(n)\mid|D|\equiv2\bmod4\}|\;+\\
\notag&&\llap{$|\{A\in\ASM(n)\mid A \text{ invariant under rotation by }\ts\frac{\pi}{2}\}|$}.\end{align}

% to fix url issues with line breaking
\let\oldurl\url
\makeatletter
\renewcommand*\url{%
        \begingroup
        \let\do\@makeother
        \dospecials
        \catcode`{1
        \catcode`}2
        \catcode`\ 10
        \url@aux
}
\newcommand*\url@aux[1]{%
        \setbox0\hbox{\oldurl{#1}}%
        \ifdim\wd0>\linewidth
                \strut
                \\
                \vbox{%
                        \hsize=\linewidth
                        \kern-\lineskip
                        \raggedright
                        \strut\oldurl{#1}%
                }%
        \else
                \hskip0pt plus\linewidth
                \penalty0
                \box0
        \fi
        \endgroup
}
\makeatother
% to fix MR issues
\gdef\MRshorten#1 #2MRend{#1}%
\gdef\MRfirsttwo#1#2{\if#1M%
MR\else MR#1#2\fi}
\def\MRfix#1{\MRshorten\MRfirsttwo#1 MRend}
\renewcommand\MR[1]{\relax\ifhmode\unskip\spacefactor3000 \space\fi
  \MRhref{\MRfix{#1}}{{\tiny \MRfix{#1}}}}
\renewcommand{\MRhref}[2]{%
 \href{http://www.ams.org/mathscinet-getitem?mr=#1}{#2}}
\bibliography{biblio}
\bibliographystyle{amsplainhyper}
\end{document}